\documentclass{article}

\usepackage{arxiv}

\usepackage[numbers, sort]{natbib}
\usepackage{url}
\usepackage{pgfplots}
\pgfplotsset{compat=1.18}
\usepackage[short]{optidef}
\usepackage{amsmath}
\usepackage{booktabs}
\usepackage{tikz}
\usepackage{amssymb}
\usepackage{stmaryrd}
\usepackage{comment}
\usepackage{bm}
\usepackage{amsfonts}
\usepackage{placeins}
\usepackage{amsthm}

\theoremstyle{plain}
\usepackage[utf8]{inputenc} 
\usepackage[T1]{fontenc}    
\usepackage{hyperref}       
\usepackage{url}            
\usepackage{booktabs}       
\usepackage{amsfonts}       
\usepackage{nicefrac}       
\usepackage{microtype}      
\usepackage{lipsum}		
\usepackage{graphicx}
\usepackage{doi}

\newcommand{\lb}[1]{\textcolor{black}{#1}}

\begin{document}

\title{A Unified Approach to Extract Interpretable Rules from Tree Ensembles via Integer Programming}


\author{ \href{https://orcid.org/0000-0002-7931-5755}{\includegraphics[scale=0.06]{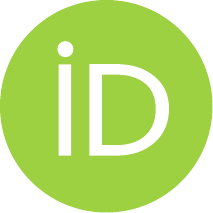}\hspace{1mm}Lorenzo Bonasera} \\
	Department of Mathematics ``Felice Casorati''\\
    University of Pavia\\
	Pavia 27100 \\
	\texttt{lorenzo.bonasera01@universitadipavia.it} \\
	\And
	\href{https://orcid.org/0000-0002-0832-8700}{\includegraphics[scale=0.06]{orcid.pdf}\hspace{1mm}Emilio Carrizosa} \\
	Institute of Mathematics of the University of Seville\\ University of Seville\\
	Seville 41012  \\
	\texttt{ecarrizosa@us.es} \\
}

\date{}

\renewcommand{\headeright}{}
\renewcommand{\undertitle}{}
\renewcommand{\shorttitle}{Extract Intepretable Rules from Tree Ensembles via Integer Programming}

\hypersetup{
pdftitle={A Unified Approach to Extract Interpretable Rules from Tree Ensembles via Integer Programming},
pdfsubject={q-bio.NC, q-bio.QM},
pdfauthor={Lorenzo Bonasera, Emilio Carrizosa},
pdfkeywords={
Machine learning,
Integer programming,
Interpretability,
Explanation fidelity,
Rule extraction,
Tree ensemble,
Time series classification},
}

\maketitle

\begin{abstract}
Tree ensembles are very popular machine learning models, known for their effectiveness in supervised classification and regression tasks. Their performance derives from aggregating predictions of multiple decision trees, which are renowned for their interpretability properties. However, tree ensemble models do not reliably exhibit interpretable output. Our work aims to extract an optimized list of rules from a trained tree ensemble, providing the user with a condensed, interpretable model that retains most of the predictive power of the full model. Our approach consists of solving a set partitioning problem formulated through Integer Programming. The proposed method works with either tabular or time series data, for both classification and regression tasks, \lb{and its flexible formulation can include any arbitrary loss or regularization functions.} Our extensive computational experiments offer statistically significant evidence that our method is competitive with other rule extraction methods \lb{in terms of predictive performance and fidelity towards the tree ensemble. Moreover, we empirically show that the proposed method effectively extracts interpretable rules from tree ensemble that are designed for time series data.}
\end{abstract}


\def\tsc#1{\csdef{#1}{\textsc{\lowercase{#1}}\xspace}}
\tsc{WGM}
\tsc{QE}
\newcommand{\parents}{Pa} 

\let\WriteBookmarks\relax
\def\floatpagepagefraction{1}
\def\textpagefraction{.001}














\section{Introduction}

The demand for interpretable and explainable machine learning models has arisen across various applications. \lb{Indeed, decision-makers need to understand the causal mechanisms behind a model's output in order to justify the choices made, especially in high-risk domains \cite{aria2021comparison}.} Consequently, research in interpretability has expanded significantly within both the machine learning and operations research communities over the past decade. Explainable Artificial Intelligence (XAI) has emerged as a rapidly growing research field~\cite{gunning2019xai}, primarily focusing on post-hoc explanation methods to facilitate the understanding of outputs provided by black-box models such as deep neural networks. Although XAI methods can be highly effective in capturing and representing the internal decision-making processes of black-box models, their application in high-risk domains is controversial~\cite{rudin2019stop}. Therefore, researchers have gradually shifted their attention to developing inherently interpretable methods.

Commonly regarded as one of the most interpretable machine learning models, decision trees have been for decades the tool of choice for providing highly intuitive and understandable predictions, thanks to their affinity with human reasoning~\cite{silva2017optimization}. However, their interpretability is counterbalanced by poor out-of-sample performance, mainly due to overfitting and inconsistent predictions. Aggregating trees into ensembles solved this problem, giving rise to machine learning models such as random forest~\cite{breiman2001random} or gradient boosted trees~\cite{friedman2001greedy, chen2015xgboost}, which significantly improved performance and robustness for both classification and regression tasks. However, combining a large number of decision trees naturally hinders their original interpretability, producing black-box models. \lb{As a consequence, research in the field of interpretable machine learning has focused on developing methods and techniques to interpret both the predictions and the internal mechanisms of tree ensemble models. More in general, methods for interpreting ensemble models such as random forests can be grouped into two categories \cite{aria2021comparison}. The first one refers to \textit{internal processing} techniques, which involve computing ensemble-specific measures to gain interpretable information about predictions. This category includes various measurement of the feature importance, such as Mean Decrease Accuracy, Mean Decrease Impurity, and Minimal Depth \cite{ishwaran2010high, seifert2019surrogate}. From a statistical inference viewpoint, feature importance consists of an estimate for the quantity that would have been obtained by making the same prediction
with unknown true response function~\cite{hooker2021unrestricted}. As discussed in previous research, model diagnostics based on feature importance can be misleading and biased towards correlated features~\cite{hooker2021unrestricted, zhou2021unbiased, strobl2008conditional}. 
The second category comprehends all the \textit{post-hoc approaches} that aim to identify and explain the relationship between tree ensembles and their outputs \cite{aria2021comparison}. These approaches include the so-called born-again tree methods, which aim to induce a single decision tree that replicates the behavior of a given tree ensemble \cite{di2024unboxing, vidal2020born, breiman1996born}. Another class of post-hoc techniques for interpreting tree ensembles consists of rule extraction methods, which involve inducing a rule-based predictor to capture the predictive power of the ensemble model \cite{di2024unboxing, aria2021comparison}.}

Arguably, rule lists (or rule sets) represent the only competitor of decision trees in terms of interpretability~\cite{yang2017scalable}, as they both consist of logical models in the form of ``if-then'' statements that are easy to understand~\cite{rudin2019stop}. A list of rules can be built from scratch or by exploiting association rules previously mined from given data~\cite{agrawal1993mining}. An interesting connection between decision trees and rule lists involves extracting rules from previously trained tree ensembles. By doing so, the search space for selecting rules is restricted to the set of branch nodes in the tree ensemble, which represent the logical conditions used to split the training data and make predictions. Thus, extracting rules from a trained tree ensemble avoids exploring impractically large search spaces. Existing methods involve solving various optimization problems to select the best performing rules, along with their corresponding optimal weights~\cite{benard2021sirus, nodeharvest, liu2023fire, friedman2008predictive}. \lb{However, these methods lack the flexibility to incorporate diverse loss functions, as they primarily focus on minimizing the squared prediction error. This limitation makes existing approaches difficult to customize and adapt to different learning tasks. Indeed, their extension to both classification and regression tasks is challenging, and their implementations are currently incompatible with complex data types such as time series, images, or text. This can become a significant drawback, especially given the growing academic and industrial interest in developing machine learning models for these types of data. For these reasons, our aim is to provide an explanation method that condenses any trained tree ensemble into an interpretable list of rules, retaining most of its predictive power and faithfully representing its inner structure, regardless of the type of data.} 

\lb{To cope with the lack of flexibility,} Mathematical Programming has recently played a fundamental role in developing interpretable and customizable methods, inducing more sparse and accurate models~\cite{ carrizosa2013supervised}. A new field belonging to operations research started focusing on formulating Mixed-Integer Linear Programming models to learn optimal classification and regression trees~\cite{bertsimas2017optimal, carrizosa2021mathematical}, as well as optimal rule lists~\cite{wang2015learning, lawless2023interpretable, wei2019generalized}. The ability to fully control how models learn while minimizing their global prediction error has made this research field growing rapidly, particularly thanks to the recent hardware and algorithmic improvements on solving Integer Programming (IP) problems~\cite{bertsimas2017optimal, lawless2023interpretable}. However, due to their complexity, most interpretable methods based on IP problems face scalability issues, both in terms of the number of records and features contained in data. In fact, dealing with large-scale data is becoming an increasingly common requirement for machine learning methods related to data-driven decision making~\cite{ullah2006knowledge}. For this reason, developing computationally efficient interpretable methods is a challenging problem. While tree ensemble and heuristic rule extraction methods do not suffer from scalability issues, they lack the expressive power and flexibility that Mathematical Programming can offer.

\subsection{Our contribution}
In light of all the above, we present as our contribution a novel method for extracting an optimal list of rules from any tree ensemble model by the means of Mathematical Programming. The proposed method relies on solving a well-known IP formulation of the set partitioning problem~\cite{balas1976set}.
The solution of such problem yields an optimal subset of rules in terms of interpretability and predictive power. The use of Mathematical Programming makes the proposed method highly flexible and easily adaptable to various applications. It can be applied to both regression and multi-class classification tasks, explicitly handling tabular and time series data. Remarkably, its formulation allows users to easily incorporate multiple custom loss functions and robustness measures. Through extensive computational experiments, we offer statistically significant evidence that our method is a valid competitor to the state-of-the-art methods for extracting rules from tree ensembles. \lb{As fidelity and interpretability are essential for producing valuable explanations that accurately describe black-box models \cite{gilpin2018exp, zhou2021evaluating}, we introduce two novel measure to assess explanation fidelity of surrogate models from tree ensembles. Obtained results show that our approach excels in terms of internal fidelity towards the ensemble model, while achieving competitive scores in external fidelity.}

\paragraph{Outline}
The remainder of this paper is organized as follows.
In Section~\ref{sec2}, we present an overview on previous works about learning rule lists from data, with a specific focus on related works on rule extraction methods from tree ensembles, discussing their main advantages and drawbacks. \lb{Additionally, we discuss the role of explanation fidelity in the context of rule extraction, and we review existing works on explainable methods for time series analysis.}
In Section~\ref{sec3}, we provide the mathematical notation and preliminary concept for introducing the problem, \lb{proposing two novel measures of internal fidelity for surrogate models from tree ensembles.}
In Section~\ref{sec4}, we introduce our framework and present the optimization problem, \lb{outlining the procedures for validating the proposed method and performing data inference.}
In Section~\ref{sec5}, we report and discuss our computational results \lb{in terms of predictive performance and explanation fidelity on classification and regression tasks using benchmark tabular datasets. Furthermore, we conduct a computational evaluation of the proposed method on time series multi-class classification tasks.}
Finally, Section~\ref{sec6} concludes the paper with a perspective on future works and further research directions.

\section{Related work} \label{sec2}

\lb{In this section, we review and discuss previous works on learning of rule sets from tabular data, with a particular focus on the past literature about rule extraction from ensemble models. Moreover, we discuss the role of explanation fidelity of surrogate models, such as rule sets, with respect to the black-box models from which they are derived. Finally, we explore related work on explainable methods for time series data, aiming to highlight the position of this paper in relation to existing research gaps.} Remarkably, earlier works about learning rule sets span multiple fields~\cite{furnkranz2012foundations}. Various approaches have been explored, including Disjunctive Normal Form learning theory \cite{valiant1984theory, klivans2001learning, feldman2012learning, su2016learning}, contrast set learning \cite{azevedo2010rules, bay1999detecting, kralj2007contrast} and frequent pattern mining \cite{dong1999efficient, garcia2010new}. The problem of learning rule lists directly from data has been well studied \cite{lakkaraju2016interpretable, domingos1996unifying}, involving Bayesian frameworks \cite{letham2015}, heuristic methods~\cite{chen2006new, cheng2006discriminative, clark1989cn2, cohen1999simple} and exact methods based on Mathematical Programming~\cite{wang2015learning, lawless2023interpretable, wei2019generalized}.

\subsection{Rule extraction from tree ensembles} Extracting rules from a trained tree ensemble avoids exploring impractically large search spaces and is considered a state-of-the-art post-hoc approach for interpreting random forests~\cite{aria2021comparison, johansson2022rule}. Friedman and Popescu~\cite{friedman2008predictive} first proposed a specialized method for extracting lists of rules from both classification and regression tree ensembles, called RuleFit. Their method extracts a subset of rules by minimizing a LASSO optimization problem, selecting the best regularized weight for each rule. However, despite its flexibility, RuleFit tends to perform poorly when extracting rules from deeper trees with highly correlated features~\cite{liu2023fire}. Meinshausen~\cite{nodeharvest} was the first to consider the tradeoff between predictive performance and interpretability with NodeHarvest. The author proposed to extract rules from a tree ensemble by solving a set partitioning problem of non-negative weights through Quadratic Programming. NodeHarvest demonstrates good predictive performance and computational efficiency, offering users an inherently interpretable and sparse model. However, its extension to classification problems is limited to handling no more than two classes. More recently, Bénard et al. proposed a rule extraction method for both classification~\cite{benard2021sirus} and regression tasks~\cite{benard2021interpretable} called SIRUS. Both variants of SIRUS involve a preprocessing step where random forests are restricted to split nodes based on the~$q$-empirical quantiles for each feature, with~$q$ being a user-defined parameter. Then, the extracted rules are filtered based on an user-defined threshold frequency~$p_0$. These filtered rules are then aggregated either by taking their average prediction or by solving a ridge regression problem, depending on the task at hand. Both variants demonstrate good predictive performance and robustness to data perturbation. \lb{Nonetheless, its classification variant cannot handle multi-class classification problems.} Despite providing interpretable output, the aforementioned methods are incapable of combining extracted rules by imposing an order on features, failing to resemble a tree-like structure. This issue has been addressed by Liu and Mazumder with FIRE~\cite{liu2023fire}, which outperforms existing methods in terms of regression performance and interpretability. The authors proposed a quadratic problem introducing a non-smooth fusion penalty. Their framework promotes solutions that exhibit a tree-like structure, thereby improving interpretability and sparsity of extracted rules. The problem is then solved through an ad-hoc optimization algorithm that leverages its block structure. However, while FIRE excels in regression performance, its extension to classification tasks is not straightforward and has not been investigated by the authors. Moreover, FIRE \lb{does not offer a direct control on the number of extracted rules.}

\lb{\subsection{Explanation fidelity} In the literature, the discrepancy between a black-box (or opaque) model and a surrogate model extracted from it is commonly referred to as explanation \textit{fidelity}~\cite{messalas2019, velmurugan2021, johansson2022rule} or \textit{completeness}~\cite{gilpin2018exp}.  A surrogate model whose predictions closely match those of the opaque model is considered highly faithful. In the context of rule extraction, explanation fidelity requires the extracted rules to faithfully represent the behavior of the trained black-box model. As argued by Zhou~\cite{zhou2004rule} and Johansson et al. \cite{johansson2022rule}, rule extraction methods should be used to accomplish two distinct tasks. The first one is to understand the underlying reasoning behind individual predictions made by a black-box model, while the second one consists of extracting an interpretable model that serves as a more understandable predictor. Methods handling the former task are evaluated based on fidelity, whereas out-of-sample accuracy is the primary criterion for methods focused on the latter task. Related works on rule extraction methods from tree ensembles fall into the second category. Indeed, the reported results are compared based on predictive performance, while any form of fidelity measure is neglected. According to Gilpin et al.~\cite{gilpin2018exp}, fidelity is necessary to produce explanations that accurately describe the black-box model. In fact, interpretability and fidelity are both needed to obtain valuable explanations~\cite{zhou2021evaluating}. However, they can conflict, as developing an accurate surrogate often requires a more complex model, which can reduce its interpretability \cite{zhou2004rule, gilpin2018exp, johansson2022rule, messalas2019}. Therefore, a tradeoff between interpretability and fidelity comes into play, and explanation methods should not be evaluated based on a single point along this tradeoff \cite{gilpin2018exp}.
As first discussed by Messalas et al. \cite{messalas2019}, two types of fidelity can be identified. The first, known as \textit{external} fidelity \cite{velmurugan2023through}, refers to the degree to which the predictions of the surrogate model agree with those of the opaque model. Common measures of external fidelity include calculating the disagreement between the two models, often expressed as the fraction of unexplained variance or through mean squared error. While external fidelity evaluates the consistency between the predictions of the surrogate and the opaque models, it does not guarantee that the surrogate’s decision-making process is faithful to the original one \cite{messalas2019}. This aspect is captured by the \textit{internal} fidelity, which is notably hard to assess \cite{velmurugan2023through}. The model-agnostic approach to assessing internal fidelity involves measuring the alignment between the feature importance rankings provided by the surrogate and opaque models. For example, this can be done by sorting the features of a dataset according to their importance as determined by both models, and then computing the rank correlation (e.g., Kendall's $\tau$) between the two orderings \cite{velmurugan2023through, messalas2019}. Currently, there is a lack of model-specific methods for assessing the internal fidelity of rule extraction techniques applied to tree ensembles.}

\subsection{Explainable methods and tree ensembles for temporal data}
Explainable machine learning methods for temporal data mainly focus on providing users with explanations derived from the data itself, such as time points, subsequences, or instance-based methods~\cite{theissler2022explainable}. Indeed, the majority of earlier works about knowledge extraction from temporal data in terms of logical rules involve transitions between time points~\cite{ullah2006knowledge, lu1998stock}, subsequences~\cite{ozden1998cyclic, das1998rule, kadous1999learning, cotofrei2002rule} or Allen's temporal logic~\cite{sciavicco2019towards, allen1983maintaining, bonasera2024learning}. Explainable methods designed to discover knowledge from time series include various extensions of the decision tree induction algorithms. For instance, Sciavicco and Stan developed a time series extension of the C4.5 algorithm~\cite{sciavicco2020knowledge}, following the work of Brunello et al.~\cite{brunello2019interval} extending the ID3 algorithm. Similarly, previous works introduced tree ensemble models specifically designed to operate with time series, such as interval-based boosted classifiers~\cite{rodriguez2001boosting} and random forests~\cite{deng2013time}, time points-based random forests~\cite{auret2010change, dudek2015short} and Proximity Forests~\cite{lucas2019proximity}. Arguably, the most explainable tree ensemble model for time series is represented by the generalized Random Shapelet Forest (gRSF), introduced by Karlsson et al.~\cite{KarlssonForest}. The authors were the first to propose a tree ensemble classification and regression method based on \textit{shapelets}~\cite{YeKeogh1, LearningSHAPLET}, which are generally defined as class-discriminative subsequences of time series. 
More specifically, these subsequences are the ones that vary the most among different classes of the dataset. The key idea of shapelets-based methods is to extract interesting features by mining local subsequences from time series datasets and utilize them for learning tasks such as classification, regression and clustering. Since shapelets are short time series subsequences, they can be visually inspected, providing users with interpretable output.
While the generalized Random Shapelet Forest (gRSF) framework is among the state-of-the-art explainable methods for temporal data, its underlying machine learning model is not interpretable due to its ensemble structure. To the best of our knowledge, our work represents the first rule set extraction method \lb{whose implementation is directly} capable of handling tree ensembles designed for temporal data.

\section{Preliminaries} \label{sec3}

In this section, we introduce the main notation and formalism about time series, decision trees and tree ensembles, to prepare the ground for presenting our rules extraction problem. Our formalism allows us to represent trees and ensemble models as \lb{feature-agnostic} sets of rules for both classification and regression settings. \lb{We also present two novel measures of internal fidelity for surrogate models based on rule extraction from tree ensembles.}

\subsection{Time series}
A time series~$\bm{x} = [x_1, x_2, \dots, x_P]$ is an ordered collection of~$P$ real values.
Given a time series~\mbox{$\bm{x} \in \mathbb{R}^P$}, a time series subsequence of~$\bm{x}$ of length~$J \leq P$ starting at position~$p$ is the ordered collection of values~\mbox{$\bm{x}_p^J = [x_p, \dots, x_{p+J-1}]$}. Please note that there are a total of~$P-J+1$ subsequences.
Given two time series~$\bm{x}, \bm{z} \in \mathbb{R}^P$, we consider as time series distance the Euclidean norm between the two corresponding vectors~$ || {\bm{x} - \bm{z}} || _2$. In a general setting, distance measures different from the Euclidean norm are allowed.
Given two time series~$\bm{x} \in \mathbb{R}^P$ and~$\bm{z} \in \mathbb{R}^J$, with~$J \leq P$, the time series \textit{subsequence distance} is the minimum distance between~$\bm{z}$ and any subsequence~$\bm{x}_p^J$ of~$\bm{x}$, that is 

\begin{equation}
\label{eq:subdist}
    \mbox{dist}(\bm{x}, {\bm{z}}) = \min_{1\leq p \leq P - J + 1} || {\bm{x}_p^J - \bm{z}} ||_2.
\end{equation}

\subsection{Decision trees}
We consider standard supervised classification and regression problems, with a given a training dataset of~$N$~i.i.d. samples \mbox{$\mathcal{D} = \{(\bm{x}_i, y_i), i = 1, \dotsc, N\}$}, where data points~$\bm{x}_i \in \mathbb{R}^P$ are expressed as vectors containing~$P$ tabular features or time series of length~$P$. In the case of classification task, targets~$y_i \in \mathcal{C}$ represent a categorical variable, where~$\mathcal{C}$ is a set of class labels. In the case of regression task, targets~$y_i \in \mathbb{R}$ represent a numerical variable. The given dataset is used to construct a predictor that maps each data point to the codomain of the corresponding task. A binary decision tree of maximal depth~$D$ is a predictor that recursively divides the feature space into at most~$2^D$ non-overlapping partitions called \textit{leaf nodes}. Each partition is defined by an ordered collection of disjunctive splittings called \textit{branch nodes}. Splittings are defined by a tuple~$(a_t(\bm{x}), b_t)$, for each~$t \in \mathcal{B}$, where~$\mathcal{B}$ is the set of branch nodes. The first element of the tuple represents the feature selector, whereas the second element corresponds to the threshold value. Each splitting compares the selected feature to the threshold value. If a data point exhibits a value lower or equal than the threshold, it is assigned to the left child of the branch node; otherwise, it is assigned to the right one. In the case of univariate decision trees for tabular data, the feature selector is the scalar product between point~$\bm{x} \in \mathbb{R}^P$ and one canonical base~$\bm{e}^{(p)}$ of the same feature space~$\mathbb{R}^P$. Figure~\ref{fig:tree} depicts an example of an univariate binary decision tree of maximal depth~$D=2$ for the classification of tabular data. It contains branch nodes~$\mathcal{B} = \{1, 2, 3\}$ and leaf nodes~$\mathcal{L} = \{1, 2, 3, 4\}$. The root node selects the 10th feature of each point contained in the dataset, and compares each value to the threshold~$b_1 =0.7$. The other branch nodes operate in a similar way. Each leaf node is associated to a colour representing one of two class labels in~$\mathcal{C} = \{\text{\textcolor{red}{Class 0}}, \text{\textcolor{blue}{Class 1}}\}$.

\begin{figure}[!htb]
\centering
\small
\begin{tikzpicture}[level distance=2cm,
level 1/.style={sibling distance=7.5cm},
level 2/.style={sibling distance=3cm}]
\node[circle,draw,fill=lightgray!20,text=black] (1) {1}
child {node[circle,draw,fill=lightgray!20,text=black] (2) {2}
child {node[circle,draw,fill=red,text=white] (4) {1}}
child {node[circle,draw,fill=blue,text=white] (5) {2}}
}
child {node[circle,draw,fill=lightgray!20,text=black] (3) {3}
child {node[circle,draw,fill=red,text=white] (6) {3}}
child {node[circle,draw,fill=blue,text=white] (7) {4}}
};
\path (1) -- (2) node [pos=0.1, align=center, left=0.4cm, font=\normalsize, color=black] {$ \langle \bm{e}^{(10)}, \bm{x} \rangle \leq 0.7$};
\path (1) -- (3) node [pos=0.1, align=center, right=0.4cm, font=\normalsize, color=black] {$ \langle \bm{e}^{(10)}, \bm{x} \rangle > 0.7$};
\path (2) -- (4) node [pos=0.1, align=center, left=0.3cm, font=\normalsize, color=black] {$ \langle \bm{e}^{(8)}, \bm{x} \rangle \leq 12.2$};
\path (2) -- (5) node [pos=0.1, align=center, right=0.3cm, font=\normalsize, color=black] {$ \langle \bm{e}^{(8)}, \bm{x} \rangle > 12.2$};
\path (3) -- (6) node [pos=0.1, align=center, left=0.3cm, font=\normalsize, color=black] {$ \langle \bm{e}^{(2)}, \bm{x} \rangle \leq 97.8$};
\path (3) -- (7) node [pos=0.1, align=center, right=0.3cm, font=\normalsize, color=black] {$\langle \bm{e}^{(2)}, \bm{x} \rangle > 97.8$};
\end{tikzpicture}
\caption{Example of a binary decision tree of depth 2 for the classification of tabular data.}
\label{fig:tree}
\end{figure}

\subsubsection{Shapelets-based trees}
\begin{figure}[htb]
\centering
\includegraphics[width=0.95\textwidth]{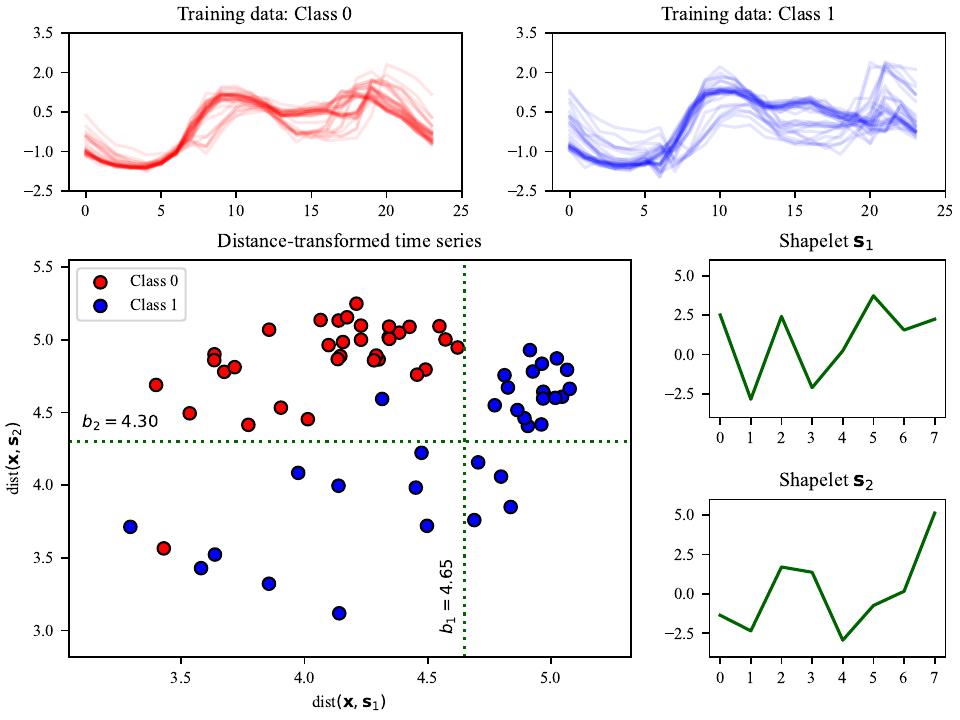}
\caption{Example of shapelets-based decision tree built upon the \textit{ItalyPowerDemand} dataset from the UCR repository~\cite{UCRArchive2018}. The upper subplots show the training time series in red (\textcolor{red}{Class 0}) and blue (\textcolor{blue}{Class 1}). The right subplots show the two shapelets~$\bm{s}_1, \bm{s}_2$ in dark green. The left subplot displays the scatter plot of the time series based on dist$(\bm{x}, \bm{s}_1)$ and dist$(\bm{x}, \bm{s}_2)$.\label{fig:dummy}}
\end{figure}

If points contained in the dataset are time series, the first element of each splitting~$(a_t(\bm{x}), b_t)$ represents an extractor of temporal features. Without loss of generality, we consider as temporal feature the distance~\eqref{eq:subdist} computed between time series and meaningful subsequences sampled from the dataset.
We show how to use data-mined shapelets as subsequences to learn a decision tree for time series interpretable classification~\cite{YeKeogh1, bonasera2024optimal}. As an example, we use the \textit{ItalyPowerDemand} dataset from the UCR repository~\cite{UCRArchive2018}. The training set includes 67 time series, each 24 units long, representing twelve monthly electrical power demand patterns from Italy. The classification task is to distinguish days between October and March from those between April and September, which correspond to the two class labels in~$\mathcal{C} = \{\text{\textcolor{red}{Class 0}}, \text{\textcolor{blue}{Class 1}}\}$. The top two subplots in Figure~\ref{fig:dummy} show the training data for the two classes. Using one of the available data mining approaches~\cite{tslearn, LearningSHAPLET}, it is possible to extract shapelets~$\bm{s}_1$ and~$\bm{s}_2$ from the training data. These shapelets are depicted in dark green in the rightmost part of Figure~\ref{fig:dummy}.
By computing distance~\eqref{eq:subdist} between each time series~$\bm{x}$ and each of the two shapelets~$\bm{s}_1$ and~$\bm{s}_2$, we can represent each time series as a point in the scatter plot in Figure~\ref{fig:dummy}. The horizontal axis represents dist$(\bm{x}, \bm{s}_1)$, while the vertical axis represents dist$(\bm{x}, \bm{s}_2)$.
Finally, using the two distances as features for each time series, we can build a decision tree that partitions the feature space associating threshold values~$b_1$ and~$b_2$ with shapelets~$\bm{s}_1$ and~$\bm{s}_2$, respectively. Figure~\ref{fig:completeshap} shows a shapelets-based decision tree of maximal depth $D = 2$, with branch nodes~$\mathcal{B} = \{1, 2\}$ and leaf nodes~$\mathcal{L} = \{1, 2, 3\}$. Threshold values correspond to those shown in the scatter plot in Figure~\ref{fig:dummy}.

\begin{figure}[htb]
\centering

\begin{tikzpicture}[level distance=2.5cm,
level 1/.style={sibling distance=4.5cm},
level 2/.style={sibling distance=3.5cm}]

\node[circle,draw,fill=lightgray!20,text=black] (1) {1}
child {node[circle,draw,fill=lightgray!20,text=black] (3) {2}
child {node[circle,draw,fill=blue,text=white] (4) {1}}
child {node[circle,draw,fill=red,text=white] (5) {2}}}
child {node[circle,draw,fill=blue,text=white] (2) {3}
};

\path (1) -- (3) node [pos=0.1, align=center, left=0.4cm, font=\normalsize, text=black] {dist$(\bm{x}, \bm{s}_1) \leq 4.65$};
\path (1) -- (2) node [pos=0.1, align=center, right=0.4cm, font=\normalsize, text=black] {dist$(\bm{x}, \bm{s}_1) > 4.65$};
\path (3) -- (4) node [pos=0.1, align=center, left=0.3cm, font=\normalsize, text=black] {dist$(\bm{x}, \bm{s}_2) \leq 4.30$};
\path (3) -- (5) node [pos=0.1, align=center, right=0.3cm, font=\normalsize, text=black] {dist$(\bm{x}, \bm{s}_2) > 4.30$};
\end{tikzpicture}

\caption{Example of a shapelets-based decision tree of depth 2 for the classification of temporal data.}
\label{fig:completeshap}
\end{figure}

\subsubsection{Rule extraction} We define as \textit{condition} the splitting of a given branch node~$t \in \mathcal{B}$ paired with the left ($\leq$) or right ($>$) sign. For instance, we indicate with~$(a_1(\bm{x}), b_1, \leq)$ the condition that a given point~$\bm{x}$ has to satisfy to get to the left child of the root node, that is whether or not the inequality~$a_1(\bm{x}) \leq b_1$ holds. Starting from the root node, each point is recursively evaluated at branch nodes until it is finally assigned to a leaf node to get a prediction. In the case of a classification tree, the prediction corresponds to the most frequent class contained in the leaf node. For regression, the prediction is the average value of the target variables of the assigned data points.
Since the path from the root node to each leaf node is unique, we can represent each leaf node~$l \in \mathcal{L}$ as a \textit{rule}, that is the ordered collection of conditions along its path
\begin{equation}
    \mathcal{R}_l = [(a_1(\bm{x}), b_1, \lozenge), \dotsc, (a_{\parents(l)}(\bm{x}), b_{\parents(l)}, \lozenge)],
\end{equation}
where~$\lozenge$ stands for the sign of each condition, while~$\parents(l) \in \mathcal{B}$ represents the parent node of leaf node~$l \in \mathcal{L}$. Since splittings are disjunctive, each point~$\bm{x}$ that satisfies all the conditions of a rule is uniquely assigned to the corresponding leaf node. Given a rule~$\mathcal{R}_l$, we define the set that contains all the corresponding splitting as follows
\begin{equation}
\label{splittings}
    \mathcal{S}_l = \{(a_1(\bm{x}), b_1), \dotsc, (a_{\parents(l)}(\bm{x}), b_{\parents(l)})\}.
\end{equation}
As example, let us consider the decision tree depicted in Figure~\ref{fig:tree}. The rule corresponding to the first leaf node is 
\[\mathcal{R}_1 = [(\langle \bm{e}^{(10)}, \bm{x} \rangle, 0.7, \leq), (\langle \bm{e}^{(8)}, \bm{x} \rangle, 12.2, \leq)],\] whereas the corresponding set of splittings is \[\mathcal{S}_1 = \{(\langle \bm{e}^{(10)}, \bm{x}, \rangle, 0.7), (\langle \bm{e}^{(8)}, \bm{x} \rangle, 12.2)\}.\]
Finally, we can represent a decision tree~$\mathcal{T}$ by all the paths it contains, expressing it as the set containing the corresponding rules
\begin{equation}
    \mathcal{T} = \{\mathcal{R}_l \mid l \in \mathcal{L}\}.
\end{equation}

For instance, the decision tree in Figure~\ref{fig:tree} can be represented as the following set of rules
\begin{align*}
    \mathcal{T} = \{ & [(\langle \bm{e}^{(10)}, \bm{x} \rangle, 0.7, \leq), (\langle \bm{e}^{(8)}, \bm{x} \rangle, 12.2, \leq)],\\
& [(\langle \bm{e}^{(10)}, \bm{x} \rangle, 0.7, \leq), (\langle \bm{e}^{(8)}, \bm{x} \rangle, 12.2, >)],\\
& [(\langle \bm{e}^{(10)}, \bm{x} \rangle, 0.7, >), (\langle \bm{e}^{(2)}, \bm{x} \rangle, 97.8, \leq)],\\
& [(\langle \bm{e}^{(10)}, \bm{x} \rangle, 0.7, >), (\langle \bm{e}^{(2)}, \bm{x} \rangle, 97.8, >)]
\}.
\end{align*}
\lb{whereas the decision tree for time series in Figure \ref{fig:completeshap} can be represented as
\begin{align*}
    \mathcal{T} = \{
& [(\text{dist}(\bm{x}, \bm{s}_1), 4.65, \leq), (\text{dist}(\bm{x}, \bm{s}_2), 4.30, \leq)],\\
& [(\text{dist}(\bm{x}, \bm{s}_1), 4.65, \leq), (\text{dist}(\bm{x}, \bm{s}_2), 4.30, >)] \\
& [(\text{dist}(\bm{x}, \bm{s}_1), 4.65, >)]
\}.
\end{align*}}

\subsection{Tree ensembles}
A tree ensemble classification or regression model consists of learning multiple decision trees by bagging~\cite{breiman1996bagging} or boosting~\cite{friedman2001greedy}, leveraging their aggregate predictions to obtain more accurate and robust results. The former approach consists of training a collection of decision trees with different samples bootstrapped from the same training dataset, obtaining a random forest~\cite{breiman2001random}. Instead, in the latter approach trees are sequentially trained by exploiting the prediction error of the previous ones, obtaining boosted learners like XGBoost~\cite{chen2015xgboost}. For sake of simplicity, in the reminder of this paper we focus on random forests as representative of tree ensemble models. Exploiting the introduced formalism, we can express a random forest~$RF$ containing decision trees~$\mathcal{T}_k$, with~$k \in K$ as their disjoint union
, that is
\begin{equation}
    RF = \bigsqcup_{k \in K} \mathcal{T}_k.
\end{equation}

In the classification setting, the prediction of a random forest corresponds to the majority vote of its decision trees. For regression tasks, the prediction is the average of all the predictions obtained by its decision trees. A generalized Random Shapelet Forest~\cite{KarlssonForest} is obtained by bagging decision trees whose feature selectors correspond to~\eqref{eq:subdist} applied to sampled shapelets. Similarly, a Time Series Forest~\cite{deng2013time} can be obtained by applying the same procedure choosing as feature one among mean, standard deviation and slope computed on sampled intervals. 
With this formalism, we can represent any random forest as a set of rules. 

\lb{\subsection{Rule fidelity} \label{sec3.4} Our goal is to extract a surrogate model from a trained random forest in the form of a set (or list) $\mathcal{Z} = \{\mathcal{R}_l\}_{l=1}^{\ell}$ containing $\ell$ rules. By doing so, we aim to effectively retain the predictive power of the tree ensemble and faithfully represent its inner structure in an interpretable way. 
To measure the internal fidelity of a set of rules towards the tree ensemble, we introduce two novel measures of representativeness. The first one is the percentage of trees whose paths are represented by a list of rules.
A tree $ \mathcal{T}$ is \textit{path-represented} by a list of rules $\mathcal{Z}$ if $\mathcal{Z} \cap \mathcal{T} \neq \emptyset$. In other words, a tree of the random forest is path-represented by a list if the latter contains at least one of the paths from the root node of the tree to its leaf nodes. The second fidelity measure corresponds to the percentage of trees whose conditions (or branch nodes) are represented by a list of rules. A tree $\mathcal{T}$ is \textit{node-represented} by a list of rules $\mathcal{Z}$ if the latter contains at least one logical condition (or branch node) of the tree. In this way, we can measure how closely a list of rules resembles the random forest, depending on the amount of path-represented and node-represented trees. For example, the tree in Figure~\ref{fig:tree}
\begin{align*}
    \mathcal{T} = \{ & [(\langle \bm{e}^{(10)}, \bm{x} \rangle, 0.7, \leq), (\langle \bm{e}^{(8)}, \bm{x} \rangle, 12.2, \leq)],\\
& [(\langle \bm{e}^{(10)}, \bm{x} \rangle, 0.7, \leq), (\langle \bm{e}^{(8)}, \bm{x} \rangle, 12.2, >)],\\
& [(\langle \bm{e}^{(10)}, \bm{x} \rangle, 0.7, >), (\langle \bm{e}^{(2)}, \bm{x} \rangle, 97.8, \leq)],\\
& [(\langle \bm{e}^{(10)}, \bm{x} \rangle, 0.7, >), (\langle \bm{e}^{(2)}, \bm{x} \rangle, 97.8, >)]
\}
\end{align*}
is path-represented and node-represented by the rule
\[\mathcal{R}_1 = [(\langle \bm{e}^{(10)}, \bm{x} \rangle, 0.7, \leq), (\langle \bm{e}^{(8)}, \bm{x} \rangle, 12.2, \leq)],\]
while it is only node-represented by the rule
\[\mathcal{R}_2 = [(\langle \bm{e}^{(10)}, \bm{x} \rangle, 0.7, \leq), (\langle \bm{e}^{(5)}, \bm{x} \rangle, -4.1, \leq)]
,\]
as they only share the condition $(\langle \bm{e}^{(10)}, \bm{x} \rangle, 0.7, \leq)$.
}

\section{Methodology} \label{sec4}

In this section, we present and discuss the proposed method. Given a trained random forest containing a total of~$L$ rules, our approach consists of selecting an optimal list of rules $\mathcal{Z}$ by solving a set partitioning problem. In this way, we can guarantee that each point is assigned to only one of the rules it satisfies, mimicking \lb{the partition of the feature space provided by} decision trees and enhancing overall interpretability. 
The key idea is to select the best subset of rules according to their robustness to data perturbation and predictive power, measured by computing rule-specific regularization and loss function, respectively.

\subsection{Preprocessing step}\label{sec:4.1}

Since the set partitioning problem is NP-hard~\cite{lewis1983michael}, we alleviate the computational burden of our approach by separating the computation of regularization and loss functions from the optimization problem. Therefore, we characterize each rule by three quantities that are computed during the following preprocessing step. We define the \textit{stability} of the~$j$-th rule as the \lb{weighted} sum of the S{\o}rensen-Dice index~\cite{chao2006abundance, li2020generic} computed over all the rules contained in the random forest, ignoring their inequality sign. \lb{Given the set of splittings in the form \eqref{splittings}, for each $l = 1, \dotsc, L$, where $L$ is the amount of rules in the random forest, we compute the stability $\phi_j$ of rule $R_j$ as follows
\begin{equation}
        \phi_j = \sum_{l \neq j, \, l = 1}^L      \lvert w_l \rvert \, \frac{2 \lvert \mathcal{S}_j \cap \mathcal{S}_l \rvert}{\lvert \mathcal{S}_j \rvert + \lvert \mathcal{S}_l \rvert},
\end{equation}
where $w_l$ is the weight of the $l$-th rule in case of weighted tree ensembles ($w_l = 1$ if unweighted).} This measures the total amount of shared splittings between rule~$\mathcal{R}_j$ and all the other rules in the tree ensemble. 
\lb{As example, consider the following rules extracted from a random forest
\begin{align*}
    & \mathcal{R}_1 = [(\langle \bm{e}^{(10)}, \bm{x} \rangle, 0.7, \leq), (\langle \bm{e}^{(8)}, \bm{x} \rangle, 12.2, \leq)],\\
& \mathcal{R}_2 = [(\langle \bm{e}^{(10)}, \bm{x} \rangle, 0.7, \leq), (\langle \bm{e}^{(8)}, \bm{x} \rangle, 12.2, >)],\\
& \mathcal{R}_3 = [(\langle \bm{e}^{(10)}, \bm{x} \rangle, 0.7, >)],
\end{align*}
with their corresponding sets of splittings
\begin{align*}
    & \mathcal{S}_1 = \{(\langle \bm{e}^{(10)}, \bm{x} \rangle, 0.7), (\langle \bm{e}^{(8)}, \bm{x} \rangle, 12.2)\},\\
& \mathcal{S}_2 = \{(\langle \bm{e}^{(10)}, \bm{x} \rangle, 0.7), (\langle \bm{e}^{(8)}, \bm{x} \rangle, 12.2)\},\\
& \mathcal{S}_3 = \{(\langle \bm{e}^{(10)}, \bm{x} \rangle, 0.7)\}.
\end{align*}
The stability of each rule is computed as follows
\begin{align*}
    & \phi_1 = \frac{2 \lvert \mathcal{S}_1 \cap \mathcal{S}_2 \rvert}{\lvert \mathcal{S}_1 \rvert + \lvert \mathcal{S}_2 \rvert} +  \frac{2 \lvert \mathcal{S}_1 \cap \mathcal{S}_3 \rvert}{\lvert \mathcal{S}_1 \rvert + \lvert \mathcal{S}_3 \rvert} = \frac{2(2)}{2 + 2} + \frac{2(1)}{2 + 1} = 1.67,\\
& \phi_2 = \frac{2 \lvert \mathcal{S}_2 \cap \mathcal{S}_1 \rvert}{\lvert \mathcal{S}_2 \rvert + \lvert \mathcal{S}_1 \rvert} +  \frac{2 \lvert \mathcal{S}_2 \cap \mathcal{S}_3 \rvert}{\lvert \mathcal{S}_2 \rvert + \lvert \mathcal{S}_3 \rvert} = \frac{2(2)}{2 + 2} + \frac{2(1)}{2 + 1} = 1.67,\\
& \phi_3 = \frac{2 \lvert \mathcal{S}_3 \cap \mathcal{S}_1 \rvert}{\lvert \mathcal{S}_3 \rvert + \lvert \mathcal{S}_1 \rvert} +  \frac{2 \lvert \mathcal{S}_3 \cap \mathcal{S}_2 \rvert}{\lvert \mathcal{S}_3 \rvert + \lvert \mathcal{S}_2 \rvert} = \frac{2(1)}{1 + 2} + \frac{2(1)}{1 + 2} = 1.33.
\end{align*}
}
In other words, we measure how much a logical condition is replicated through the whole ensemble. A rule made of conditions that are consistently replicated inside the random forest is more stable than a rule containing less frequent conditions. By maximizing their stability, selected rules are less sensitive to data perturbation, and they can more likely resemble a decision tree structure. This approach combines the effects of the fusion penalty of FIRE~\cite{liu2023fire} and the occurrence frequency of SIRUS~\cite{benard2021interpretable}. Indeed, the former penalty term promotes the learning of rules with more shared splittings, improving their interpretability. Instead, the latter method prunes the rules that entirely occur less than the frequency threshold parameter~$p_0$.
In order to quantify the predictive power of the~$j$-th rule, we compute its loss depending on the learning task as follows
\begin{equation}
    \xi_j = \begin{cases} 
        N_j - \max_{c \in \mathcal{C}} N_{cj}, & \text{if classification,}\\
        \text{MSE}(\mathcal{R}_j), & \text{if regression,}
    \end{cases}
\end{equation}
where~$N_j$ is the number of data points that satisfy the rule, $N_{cj}$ is the number of data points of class~$c \in \mathcal{C}$ among them, and \text{MSE}$(\mathcal{R}_j)$ corresponds to the Mean Squared Error (MSE) computed on the response values~$y_i \in \mathbb{R}$ of data points that satisfy rule~$\mathcal{R}_j$.
In other words, we quantify the loss of the~$j$-th rule as the number of misclassified points in the case of classification tasks. Instead, we consider the MSE computed on assigned points as the loss of the~$j$-th rule in the case of regression tasks. By computing loss and regularization values during the preprocessing step, we manage to easily handle non-linear and multiple functions. After their computation, vectors~$\phi_j$ and~$\xi_j$ are normalized to the same scale. Lastly, we keep track of the assignment for each of the~$N$ data point to each of the~$L$ rules through the matrix~$\bm{A} \in \{0, 1\}^{N \times L}$, whose entries are computed as follows
\begin{equation}
    A_{ij} = \begin{cases}
        1, & \text{if point~$\bm{x}_i$ satisfies rule~$\mathcal{R}_j$},\\
        0, & \text{otherwise}.
    \end{cases}
\end{equation}

\subsection{Optimization problem}

After employing the preprocessing step, we formulate the problem of selecting the best list of rules as a set partitioning problem through Integer Programming. We introduce a binary decision variable for each rule of the tree ensemble, that is~$z_j \in \{0, 1\}$, with~$j = 1, \dotsc, L$. Then, we exploit the corresponding quantities~$\phi_j, \xi_j$ and assignment~$\bm{A}$ to write the following IP problem
\begin{maxi!}[2]
{z}{\lambda \sum_{j=1}^L \phi_j z_j - (1 - \lambda) \sum_{j=1}^L \xi_j z_j\label{freq:obj}}{\label{model:freq}}{}
\addConstraint{\sum_{j=1}^L A_{ij} z_j}{\label{freq:con1} = 1,\quad}{i = 1, \dotsc, N}
\addConstraint{\sum_{j=1}^L z_j}{\leq \ell\label{freq:con3}}{}
\addConstraint{z_j}{\in \{0, 1\},\quad \label{freq:con4}}{j = 1, \dotsc, L,}
\end{maxi!}
\lb{where $\lambda \in [0, 1]$ represents a weight parameter to balance loss and stability, and~$\ell \in \mathbb{N}$} represents the maximum number of rules to extract. Solving this problem produces an optimal list of rules in the form \lb{$\mathcal{Z} = \{\mathcal{R}_j \mid z_j = 1, \, j = 1, \dotsc, L\}$, containing the rules whose corresponding binary variable are included in the solution.} The objective function~\eqref{freq:obj} maximizes the stability of rules while minimizing their aggregated loss. Constraint~\eqref{freq:con1} ensures that each sample is assigned to exactly one rule.
Constraint~\eqref{freq:con3} imposes an upper bound on the number of selected rules. Finally, constraint~\eqref{freq:con4} specifies the domain of the decision variables. We remark that the proposed method is independent of how stability and loss are measured. Indeed, our approach is flexible and easy to customize, allowing the end-user to compute multiple custom, non-linear stability and loss functions during the preprocessing step described in Section~\ref{sec:4.1}. \lb{In this way, the optimization problem is agnostic to the learning task, efficiently handling both regression and multi-class classification.}

\subsection{Validation procedure}\label{sec:training}
\begin{figure}[tb]
	\centering
	\includegraphics[width=.8\textwidth]{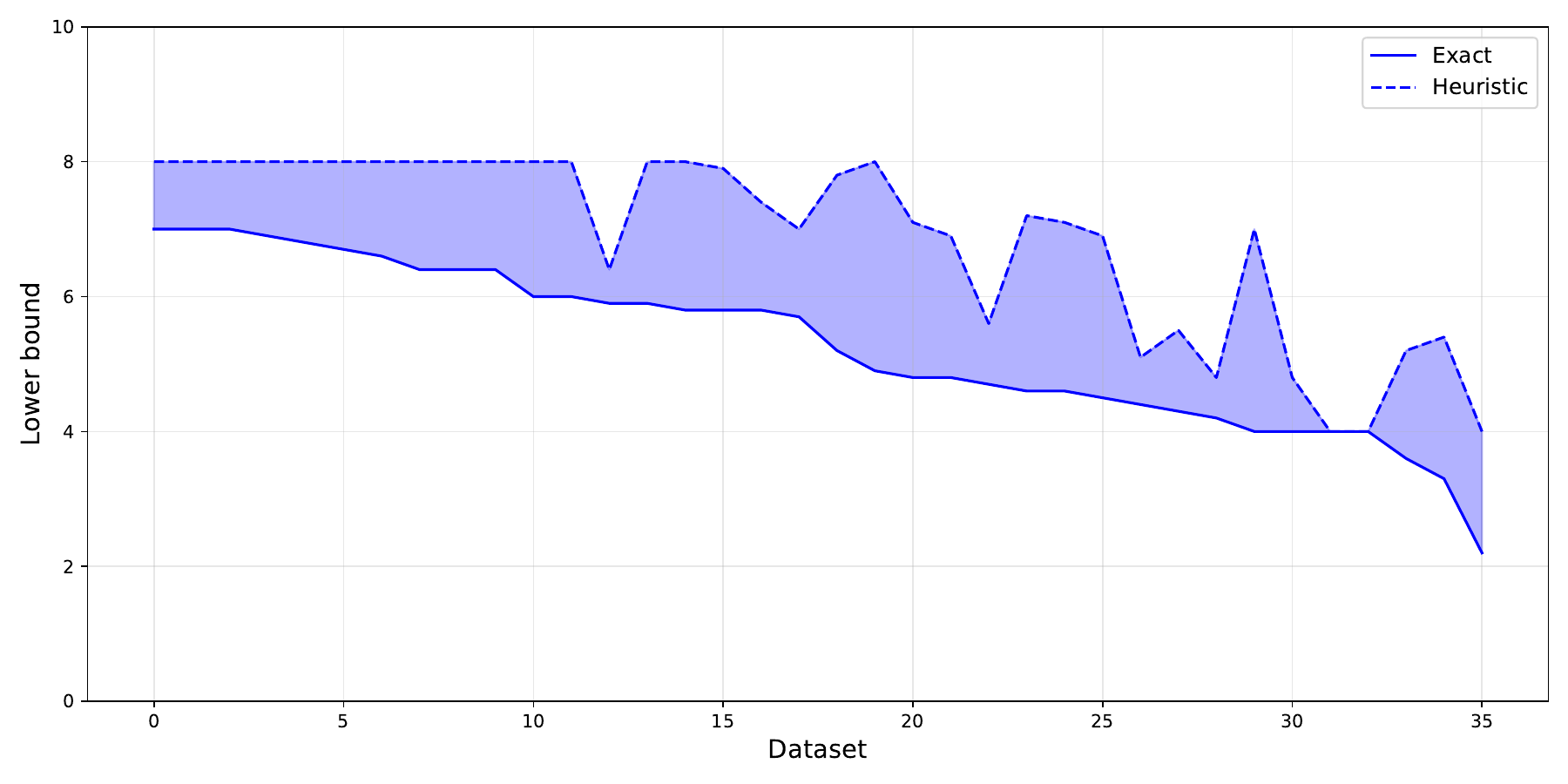}
	\caption{Comparison between exact and heuristic values for the lower bound $\underline{\ell}$ over 36 benchmark datasets.}
	\label{fig:lowerbound}
\end{figure}

In this section, we describe the procedure to train and validate the proposed method when the maximum number of rules to extract is not specified or cannot be inferred. Starting from an already trained random forest with given training and validation sets, we need to choose a proper value for parameter~$\ell$.
\lb{Therefore, we validate this parameter through an exhaustive search over a discrete set $\left \{ \underline{\ell}, \dotsc, \overline{\ell} \right \}$,} providing both a lower and an upper bound for its range. We compute the lower bound~$\underline{\ell}$ by solving the optimization problem obtained from~\eqref{model:freq} by dropping constraint~\eqref{freq:con3} and modifying the objective function as follows
\begin{mini!}[2]
{z}{\sum_{j=1}^L z_j \label{lb:obj}}{\label{model:lb}}{\underline{\ell} \; = \;}
\addConstraint{\sum_{j=1}^L A_{ij} z_j}{= 1,\label{lb:con1}\quad}{i = 1, \ldots, N}
\addConstraint{z_j}{\in \{0, 1\},\label{lb:con2}\quad}{j = 1, \ldots, L.}
\end{mini!}
We compute the upper bound~$\overline{\ell}$ by solving the same problem with the opposite objective function. The objective values obtained from these solutions correspond to the desired bounds. \lb{For possibly impractical problems, we propose two heuristic strategies to compute sufficiently tight bounds for both values. Given an ensemble of $K$ trees, we can compute a heuristic value for the lower bound~$\underline{\ell}$ as follows
\begin{equation}
\underline{\ell}^{\text{heur}} = \min_{k \in K} \, \lvert \mathcal{T}_k \rvert,
\end{equation}
that is the number of rules (or leaf nodes) contained in the smallest tree of the ensemble. This approach yields lower bounds that are close to the exact values, as shown in Figure \ref{fig:lowerbound}. To compute a heuristic value for the upper bound~$\overline{\ell}$, we can construct a decision tree with low cost complexity on the input data and count the number of generated leaf nodes. We note that obtained results can significantly vary depending on the dataset.}
After computing bounds, we solve Problem~\eqref{model:freq} for each value of~$\ell \in \left \{ \underline{\ell}, \dotsc, \overline{\ell} \right \}$ over the training set. We remark that reoptimizing the same IP problem with a different right-hand side is computationally efficient, since modern solvers can exploit effective warm start techniques~\cite{witzig2014reoptimization}. Finally, we select the value of $\ell$ that achieves the lowest loss on the validation set.
\lb{\subsection{Inference procedure}
One of the major differences between the proposed method and existing works on rule extraction consists of the lack of weighted rules. Indeed, our approach generates a list of rules that corresponds to a set partition, without incorporating any weighting strategy. We argue that this approach is more interpretable than competing ones, as it resembles a decision tree and it partitions the training data into a rule-based clustering structure, which is considered a challenging problem in the field of interpretable machine learning \cite{carrizosa2023clustering}. However, the application of an unweighted list of rules may fail on unseen data points, especially if they significantly differ from the training data. This can result in out-of-sample data falling into multiple rules or none at all, preventing the proposed method from generating a meaningful prediction. We remark that, in rule extraction for black-box models, the purpose is to explain the predictions obtained from a black-box model \cite{zhou2004rule}. Therefore, if the surrogate model fails to obtain an individual prediction as in the case of multiple or absent rules, the unseen data point can still be predicted by the tree ensemble \cite{johansson2022rule}. Nonetheless, we propose a simple heuristic to address the issue. If an unseen data point falls within two or more rules, we apply the rule that covers the largest amount of training data. Instead, if an unseen data point does not fall within any rules of the list, the prediction is the average of the response vector in case of regression, or the majority class in case of classification.}

\section{Experiments and results} \label{sec5}

In this section, we present and discuss the experimental evaluation conducted to assess the performance of the proposed method against the state-of-the-art, for regression and classification tasks of tabular data and time series classification.
To ensure a fair comparison, we selected as competing approaches the state-of-the-art methods that are reproducible. \lb{For the same reason, we excluded FIRE from the comparison, since it does not provide a fixed upper limit on the number of extracted rules. We discuss more in detail in \ref{appendix:A} the challenges of establishing a fair comparison for rule extraction methods.} We implemented all methods in Python 3.9, except for SIRUS, which was implemented in R 4.3.1. We ran all the experiments on the same machine equipped with a 3.3~GHz 8-core processor and 16 GB of memory. The proposed IP problems \eqref{model:freq} and \eqref{model:lb} are solved through Gurobi 11.0~\cite{gurobi}. \lb{From preliminary results, the impact of validating the weight parameter $\lambda$ resulted in negligible effects. Thus, in the following experimental evaluation we fixed $\lambda = 0.5$, ensuring a balanced tradeoff between loss and stability of extracted rules.} In order to compare multiple methods in terms of predictive power over the selected datasets, we employed the non-parametric Friedman test~\cite{demvsar2006statistical} combined with the signed-rank Wilcoxon-Holm post-hoc test~\cite{wilcoxon2019}, with a level of statistical significance equal to~$0.05$. 
According to Benavoli et al.~\cite{benavoli2016should}, the Wilcoxon signed-rank test is more appropriate than the more common post-hoc Nemenyi test~\cite{nemenyi}.
Details about benchmark datasets are contained in Table \ref{tab:regression}, Table \ref{tab:classtabular} and Table \ref{tab:timeseries}. During data preprocessing, we applied one-hot encoding to tabular data with categorical features and discarded samples containing missing values. For regression tasks, we standardized the response vector. No data preprocessing was applied to time series datasets.

\subsection{Tabular regression}
\begin{figure}[hb]
	\centering
	\includegraphics[width=\textwidth]{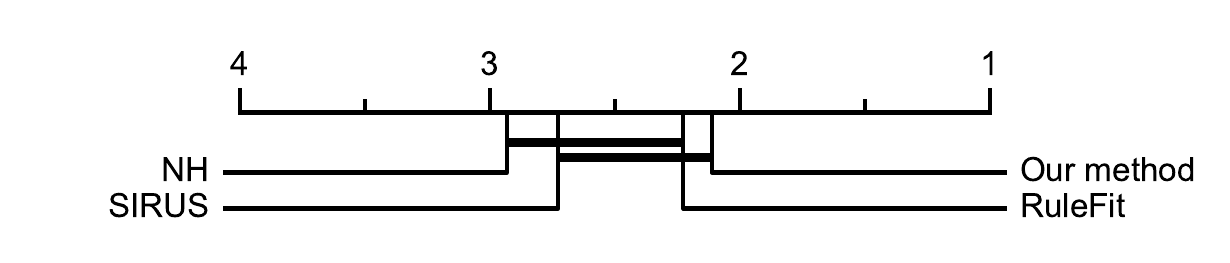}
	\caption{Critical distance diagram obtained by running the Friedman test~\cite{demvsar2006statistical} combined with the signed-rank Wilcoxon-Holm post-hoc test ~\cite{wilcoxon2019} on the results about regression of tabular data of Table~\ref{tab:regression}.}
	\label{fig:cd-regr}
\end{figure}

In this experiment, we compare the proposed method to the state-of-the-art about rule extraction for regression tasks of tabular data. In particular, competing methods are RuleFit~\cite{friedman2008predictive}, SIRUS~\cite{benard2021interpretable} and NodeHarvest (NH)~\cite{nodeharvest}. 
The experimental plan involves running each method on a selected benchmark of~22~datasets from the OpenML public repository~\cite{openml}. 
Following the experimental plan of earlier works~\cite{liu2023fire, nodeharvest}, we extract rules from a regression forest of 500 trees of depth 3 implemented through the {\tt scikit-learn} package with default settings~\cite{scikit-learn}, and fix the maximum amount of selected rules to 15. Since this quantity coincides with parameter~$\ell$, the proposed \lb{validation} procedure is unnecessary. \lb{Instead, we discard rules that cover fewer than a minimum percentage $N_{\min}$ of the training samples, where $N_{\min} \in [0.001, 0.01]$ is tuned through cross-validation. For competing methods, we fix the number of extracted rules (or nodes) to 15 and the maximal interaction order (or length of rules) to 3. Consequently, NodeHarvest does not need further parameters tuning \cite{nodeharvest}. For SIRUS, we tune the frequency threshold $p_0$ through cross-validation, and we set the $q$ parameter to 10, following the authors' guidelines \cite{benard2021interpretable}. For RuleFit, we tune the regularization parameter for the LASSO regression through cross-validation \cite{friedman2008predictive}.}
We divide each dataset into a 75/25 train-test split. We run each regression task with 30 different random seeds, and we report average and standard deviation of the out-of-sample Mean Squared Error (MSE) in Table~\ref{tab:regression}. 
\lb{From obtained results, we find that differences in average predictive performance are not statistically significant ($p \approx 0.08$). We show the pairwise statistical comparison through the critical difference diagram comparing ranks in Figure \ref{fig:cd-regr}. Remarkably, our approach achieves statistically equivalent performance to RuleFit and SIRUS, while presenting a significant improvement over NodeHarvest in terms of MSE ($p \approx 0.002$).}
The median value of computing times of the proposed method is~$\approx 7$ minutes, and the average runtime of the Gurobi IP solver is~$\approx 5\%$ of the total computing time. On the selected benchmark, competing methods are more efficient than ours. 

\begin{table}[htbp]
    \centering
        \caption{Details of datasets from the OpenML repository for tabular data regression.}
            \normalsize
    \begin{tabular}{lrr}
    \toprule
        Dataset & Records & Features\\
        \midrule
        Mercedes\_Benz & 4209 & 563  \\
        Moneyball & 114 & 54  \\ 
        abalone & 4177 & 10  \\ 
        autoMpg & 392 & 25  \\ 
        bank32nh & 8192 & 32  \\ 
        boston & 506 & 22  \\ 
        cpu\_small & 8192 & 12  \\ 
        elevators & 16599 & 18  \\ 
        house\_16H & 22784 & 16  \\ 
        houses & 20640 & 8  \\ 
        kin8nm & 8192 & 8  \\ 
        mtp & 4450 & 202  \\ 
        no2 & 500 & 7  \\ 
        pol & 15000 & 48  \\ 
        puma32H & 8192 & 32  \\ 
        socmob & 1156 & 39  \\ 
        space\_ga & 3107 & 6  \\ 
        stock & 950 & 9  \\ 
        tecator & 240 & 124  \\ 
        us\_crime & 123 & 249  \\
        wind & 6574 & 14  \\ 
        wine\_quality & 6497 & 11  \\ \bottomrule
    \end{tabular}
    \label{regression}
\end{table}

\begin{table}
\vspace{-1cm}
\centering
  \caption{Out-of-sample MSE values computed on benchmark datasets for tabular regression.}
  \label{tab:regression}
  \centering
  \small
  \begin{tabular}{l@{\hspace{1em}}*{4}{c}}
    \toprule
    Dataset & Our method & RuleFit & SIRUS & NH\\
    \midrule
    Mercedes\_Benz & 0.44~$\pm$ 0.10 & 0.46~$\pm$ 0.11 & 0.51~$\pm$ 0.10 & 0.45~$\pm$ 0.10\\
    Moneyball & 0.26~$\pm$ 0.07 & 0.22~$\pm$ 0.08 & 0.20~$\pm$ 0.07 & 0.26~$\pm$ 0.08\\
    abalone & 0.61~$\pm$ 0.06 & 0.61~$\pm$ 0.06 & 0.68~$\pm$ 0.06 & 0.60~$\pm$ 0.06\\
    autoMpg & 0.31~$\pm$ 0.06 & 0.30~$\pm$ 0.08 & 0.28~$\pm$ 0.05 & 0.32~$\pm$ 0.08\\
    bank32nh & 0.65~$\pm$ 0.03 & 0.64~$\pm$ 0.04 & 0.61~$\pm$ 0.04 & 0.65~$\pm$ 0.04\\
    boston & 0.28~$\pm$ 0.04 & 0.27~$\pm$ 0.06 & 0.29~$\pm$ 0.09 & 0.31~$\pm$ 0.06\\
    cpu\_small & 0.08~$\pm$ 0.01 & 0.08~$\pm$ 0.01 & 0.11~$\pm$ 0.01 & 0.08~$\pm$ 0.01\\
    elevators & 0.52~$\pm$ 0.03 & 0.54~$\pm$ 0.05 & 0.64~$\pm$ 0.04 & 0.53~$\pm$ 0.03\\
    house\_16H & 0.68~$\pm$ 0.05 & 0.69~$\pm$ 0.06 & 0.66~$\pm$ 0.06 & 0.69~$\pm$ 0.05\\
    houses & 0.51~$\pm$ 0.01 & 0.52~$\pm$ 0.03 & 0.56~$\pm$ 0.02 & 0.51~$\pm$ 0.01\\
    kin8nm & 0.67~$\pm$ 0.02 & 0.68~$\pm$ 0.02 & 0.65~$\pm$ 0.03 & 0.68~$\pm$ 0.05\\
    mtp & 0.74~$\pm$ 0.03 & 0.73~$\pm$ 0.04 & 0.74~$\pm$ 0.03 & 0.75~$\pm$ 0.06\\
    no2 & 0.59~$\pm$ 0.09 & 0.57~$\pm$ 0.12 & 0.54~$\pm$ 0.10 & 0.60~$\pm$ 0.13\\
    pol & 0.29~$\pm$ 0.02 & 0.28~$\pm$ 0.05 & 0.28~$\pm$ 0.02 & 0.30~$\pm$ 0.01\\
    puma32H & 0.41~$\pm$ 0.02 & 0.42~$\pm$ 0.08 & 0.56~$\pm$ 0.06 & 0.42~$\pm$ 0.01\\
    socmob & 0.33~$\pm$ 0.09 & 0.32~$\pm$ 0.11 & 0.37~$\pm$ 0.12 & 0.37~$\pm$ 0.25\\
    space\_ga & 0.53~$\pm$ 0.08 & 0.54~$\pm$ 0.10 & 0.53~$\pm$ 0.10 & 0.56~$\pm$ 0.09\\
    stock & 0.09~$\pm$ 0.01 & 0.12~$\pm$ 0.03 & 0.24~$\pm$ 0.08 & 0.10~$\pm$ 0.01\\
    tecator & 0.05~$\pm$ 0.01 & 0.06~$\pm$ 0.02 & 0.08~$\pm$ 0.02 & 0.04~$\pm$ 0.02\\
    us\_crime & 0.46~$\pm$ 0.07 & 0.33~$\pm$ 0.04 & 0.35~$\pm$ 0.09 & 0.50~$\pm$ 0.09\\
    wind & 0.36~$\pm$ 0.01 & 0.36~$\pm$ 0.03 & 0.39~$\pm$ 0.02 & 0.36~$\pm$ 0.03\\
    wine\_quality & 0.74~$\pm$ 0.04 & 0.74~$\pm$ 0.04 & 0.78~$\pm$ 0.04 & 0.75~$\pm$ 0.02\\
    \midrule
    \textit{Average} & 0.44~$\pm$ 0.04 & 0.43~$\pm$ 0.06 & 0.46~$\pm$ 0.06 & 0.45~$\pm$ 0.05\\
    \textit{Rank ($p \approx 0.08$)} & 2.11 & 2.23 & 2.73 & 2.93\\
    \bottomrule
  \end{tabular}
\end{table}
\normalsize

\lb{\subsubsection{Comparison in terms of fidelity} To assess internal fidelity, we compute the two representativeness measures introduced in Section \ref{sec3.4}, and the F1-score between the features selected by each rule extraction methods and the top $5\%$ of features ranked by their importance in the random forest, following the work of Velmurugan et al. \cite{velmurugan2023through}. To evaluate external fidelity, we report the disagreement in terms of MSE computed between the predictions of each competing method and those of the random forest. We report in Table \ref{tab:fidelity-regr} the average of each fidelity measure computed over the benchmark datasets, for each rule extraction method. For the sake of completion, we also report the average number of extracted rules. Notably, the proposed method obtains the best scores of representativeness and disagreement, while obtaining the second best F1-score. Instead, RuleFit achieves the lowest disagreement value and the best F1-score, and obtains the second best values in terms of representativeness. On average, the proposed method and NodeHarvest extract the lowest amount of rules from the tree ensemble. By contrast, SIRUS extracts the largest number of rules, due to its tendency to select rules with a very limited amount of logical conditions. We conclude that the choice of the interpretable method depends on the desired result. Indeed, SIRUS is the preferred method to obtain a list containing a large amount of short rules. On the contrary, the proposed approach is the most suitable to extract a list that contains a smaller amount of more representative rules. We show an example of extracted rules in \ref{appendix:B}.}

\begin{table}[htb]
  \caption{Average fidelity values and number of extracted rules computed on the tabular regression benchmark.}
  \label{tab:fidelity-regr}
  \small
  \centering
  \begin{tabular}{l@{\hspace{1em}}*{4}{c}}
    \toprule
    Measure & Our method & RuleFit & SIRUS & NH \\
    \midrule
    Represented trees & \textbf{0.731} & 0.607 & 0.307 & 0.514 \\
    Represented paths & \textbf{0.069} & 0.032 & 0.005 & 0.042 \\
    F1-score & 0.543 & \textbf{0.586} & 0.447 & 0.505 \\
    Disagreement & \textbf{0.047} & \textbf{0.047} & 0.092 & 0.065 \\
    \midrule
    Extracted rules & 8.44 & 10.21 & 11.50 & 7.99\\
    \bottomrule
  \end{tabular}
\end{table}

\subsection{Tabular classification}
\begin{figure}[htb]
	\centering
	\includegraphics[width=\textwidth]{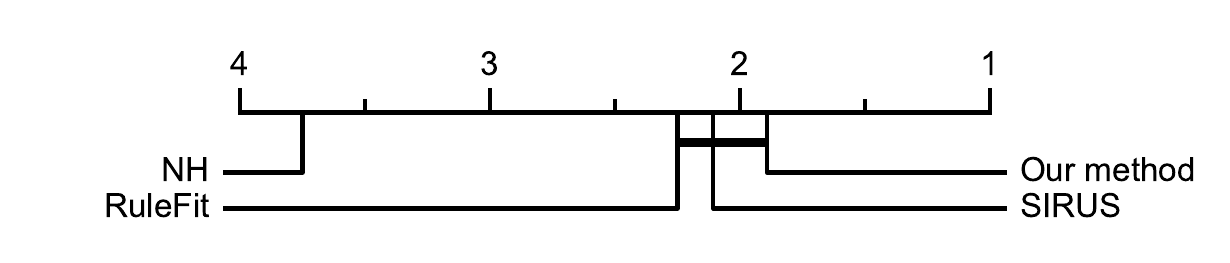}
	\caption{Critical distance diagram obtained by running the Friedman test~\cite{demvsar2006statistical} combined with the signed-rank Wilcoxon-Holm post-hoc test ~\cite{wilcoxon2019} on the results about classification of tabular data of Table~\ref{tab:classtabular}.}
	\label{fig:cd-clf}
\end{figure}

Similarly to the previous experiment, we compare the proposed method to the state-of-the-art about rule extraction for classification tasks of tabular data. As competing methods, we consider the previously mentioned approaches used for regression tasks. \lb{However, NodeHarvest and SIRUS can only handle binary classification tasks. Therefore, we select datasets with no more than~$2$ different classes.} In particular, the experimental plan involves running each method on a selected benchmark of~14~datasets from the OpenML public repository~\cite{openml}, following the approach of earlier works~\cite{lawless2023interpretable, wei2019generalized}.
To assess the effectiveness of methods in extracting very short and sparse rules that resemble a decision tree, we extract rules from a random forest of 500 trees of depth 2, and fix the maximum amount of selected rules to 4. \lb{Similar to the previous experiment, this quantity coincides with parameter~$\ell$, making the proposed validation procedure unnecessary.} \lb{For competing methods, we fix the number of extracted rules to 4 and the maximal interaction order to 2. NodeHarvest does not need further parameter tuning. For SIRUS, we tune the frequency threshold $p_0$ through cross-validation, and we set the $q$ parameter to 10. For RuleFit, we tune the regularization parameter for the logistic model through cross-validation.}
We divide each dataset into a 75/25 train-test split. We run each classification task with 30 different random seeds, reporting the average and standard deviation of the out-of-sample accuracy in Table~\ref{tab:classtabular}. \lb{From the critical difference diagram comparing ranks in Figure \ref{fig:cd-clf}, we observe that all methods achieve statistically equivalent classification performance, with the exception of NodeHarvest, which performs poorly when running with a very limited number of nodes.} The median value of computing times of the proposed method is~$\approx 1$ minute, with Gurobi's average runtime being~\mbox{$\approx 8.5\%$} of the total. Compared to competing methods, it is less efficient, as it requires more computational resources. However, our approach is more flexible than competing methods, as it can handle different \lb{loss functions and multi-class classification tasks.} 

\begin{table}[htbp]
    \centering
        \caption{Details of datasets from the OpenML repository for tabular data classification.}
            \normalsize
    \begin{tabular}{lrrr}
    \toprule
        Dataset & Records & Features & Classes\\
        \midrule
        adult & 45222 & 120 & 2  \\
        bank-marketing & 45211 & 51 & 2  \\ 
        banknote & 1372 & 4 & 2  \\ 
        blood-transfusion & 748 & 4 & 2  \\ 
        diabetes & 768 & 8 & 2  \\ 
        FICO & 9871 & 37 & 2  \\ 
        heart-c & 296 & 25 & 2  \\ 
        ilpd & 583 & 11 & 2  \\ 
        ionosphere & 351 & 34 & 2  \\ 
        MagicTelescope & 19020 & 10 & 2  \\ 
        mushroom & 5644 & 98 & 2  \\ 
        musk & 6598 & 268 & 2  \\ 
        tic-tac-toe & 958 & 27 & 2  \\ 
        wdbc & 569 & 30 & 2  \\ \bottomrule
    \end{tabular}
    \label{classification}
\end{table}

\begin{table}[!htb]
  \caption{Out-of-sample accuracy values computed on benchmark datasets for tabular classification.}
  \label{tab:classtabular}
  \small
  \centering
  \begin{tabular}{l@{\hspace{1em}}*{4}{c}}
    \toprule
    Dataset & Our method & RuleFit & SIRUS & NH \\
    \midrule
    adult & 0.81~$\pm$ 0.01 & 0.77~$\pm$ 0.02 & 0.76~$\pm$ 0.02 & 0.75~$\pm$ 0.00\\
    bank-mark. & 0.89~$\pm$ 0.00 & 0.89~$\pm$ 0.01 & 0.88~$\pm$ 0.00 & 0.88~$\pm$ 0.00\\
    banknote & 0.91~$\pm$ 0.01 & 0.85~$\pm$ 0.08 & 0.87~$\pm$ 0.03 & 0.56~$\pm$ 0.02\\
    blood-transf. & 0.74~$\pm$ 0.02 & 0.76~$\pm$ 0.02 & 0.76~$\pm$ 0.02 & 0.75~$\pm$ 0.03\\
    diabetes & 0.70~$\pm$ 0.03 & 0.70~$\pm$ 0.05 & 0.70~$\pm$ 0.04 & 0.65~$\pm$ 0.02\\
    FICO & 0.69~$\pm$ 0.01 & 0.67~$\pm$ 0.03 & 0.71~$\pm$ 0.01 & 0.52~$\pm$ 0.01\\
    heart-c & 0.73~$\pm$ 0.04 & 0.75~$\pm$ 0.06 & 0.77~$\pm$ 0.05 & 0.54~$\pm$ 0.07\\
    ilpd & 0.70~$\pm$ 0.03 & 0.72~$\pm$ 0.03 & 0.72~$\pm$ 0.03 & 0.72~$\pm$ 0.03\\
    ionosphere & 0.89~$\pm$ 0.02 & 0.83~$\pm$ 0.07 & 0.88~$\pm$ 0.04 & 0.70~$\pm$ 0.09\\
    MagicTeles. & 0.77~$\pm$ 0.02 & 0.74~$\pm$ 0.03 & 0.73~$\pm$ 0.01 & 0.65~$\pm$ 0.01\\
    mushroom & 0.91~$\pm$ 0.02 & 0.88~$\pm$ 0.07 & 0.90~$\pm$ 0.01 & 0.62~$\pm$ 0.01\\
    musk & 0.86~$\pm$ 0.01 & 0.87~$\pm$ 0.02 & 0.93~$\pm$ 0.04 & 0.85~$\pm$ 0.01\\
    tic-tac-toe & 0.68~$\pm$ 0.02 & 0.67~$\pm$ 0.05 & 0.67~$\pm$ 0.03 & 0.65~$\pm$ 0.02\\
    wdbc & 0.95~$\pm$ 0.03 & 0.92~$\pm$ 0.03 & 0.91~$\pm$ 0.02 & 0.63~$\pm$ 0.03\\
    \midrule
    \textit{Average} & 0.80~$\pm$ 0.02 & 0.79~$\pm$ 0.04 & 0.80~$\pm$ 0.03 & 0.68~$\pm$ 0.03\\
    \textit{Rank ($p < 10^{-3}$)} & \textbf{1.89} & 2.25 & 2.11 & 3.75\\ 
    \bottomrule
  \end{tabular}
\end{table}

\lb{\subsubsection{Comparison in terms of fidelity} To assess internal fidelity, we compute the same measures as in the previous experiment. For external fidelity, we calculate the disagreement as the fraction of misaligned classifications between each competing method and the random forest. We report in Table \ref{tab:fidelity-clf} the average of each fidelity measure computed over the benchmark datasets, for each rule extraction method. For the sake of completion, we also report the average number of extracted rules. From obtained results, it turns out that the proposed method excels in representativeness, while obtaining the second best F1-score. Instead, SIRUS achieves the lowest disagreement value and the best F1-score, but performs poorly in terms of representativeness. This behavior is ascribable to its tendency to extract rules with a limited number of logical conditions, as shown in \ref{appendix:B}.}

\begin{table}[htb]
  \caption{Average fidelity values and number of extracted rules computed on the tabular classification benchmark.}
  \label{tab:fidelity-clf}
  \small
  \centering
  \begin{tabular}{l@{\hspace{1em}}*{4}{c}}
    \toprule
    Measure & Our method & RuleFit & SIRUS & NH\\
    \midrule
    Represented trees & \textbf{0.355} & 0.178 & 0.064 & 0.112\\
    Represented paths & \textbf{0.021} & 0.008 & 0.006 & 0.005\\
    F1-score & 0.546 & 0.259 & \textbf{0.571} & 0.152\\
    Disagreement & 0.107 & 0.097 & \textbf{0.054} & 0.179\\
    \midrule
    Extracted rules & 3.98 & 3.24 & 3.61 & 2.48\\
    \bottomrule
  \end{tabular}
\end{table}


\subsection{Time series classification}
The aim of this experiment is to assess the effectiveness of the proposed method in handling tree ensembles specifically designed for time series. We choose the generalized Random Shapelet Forest (gRSF)~\cite{KarlssonForest} as the tree ensemble model for this purpose, providing interpretable features (shapelets).
\lb{Currently, our approach is the only rule extraction method whose implementation works on temporal data.} The experimental plan consists of running gRSF on a selected benchmark of 27 time series datasets from the UCR public repository~\cite{UCRArchive2018} following the methodology of Karlsson et al.~\cite{KarlssonForest}, and applying the proposed method to extract an optimal list of rules from it. We remark that the UCR repository provides fixed train-test splits for each dataset, avoiding the need for additional splitting strategies. For this reason, we conduct each classification task using only 10 different random seeds. The gRSF is trained using the tuned hyperparameters provided by its authors, \lb{which control the size and number of sampled shapelets. To maintain consistency with the previous experiments, we fix the maximal depth to 3.} Therefore, we select time series datasets with no more than~$2^3$ different classes. \lb{As the number of rules to extract is not fixed, we train our method by running the validation procedure with exact bounds discussed in Section~\ref{sec:training}, and we apply 5-fold cross-validation.} The average and standard deviation of the out-of-sample accuracy are reported in Table~\ref{tab:timeseries}. \lb{On average, results show that the gRSF (Full model) outperforms the proposed method, while the latter manages to extract an average of 6 interpretable rules.} The median computation time with the training procedure is~\mbox{$\approx 17$} minutes. On average, the runtime of the Gurobi IP solver is negligible. Remarkably, with given hyperparameters, gRSF is significantly faster, running in less than 30 seconds per task. Nevertheless, the computational burden of our method is justified by its interpretability and condensed predictive power. \lb{In \ref{appendix:B}, we show an example of a shapelet-based list of rules extracted from gRSF.}

\begin{table}[htb]
    \centering
        \caption{Details of datasets from the UCR repository for time series classification.}
    \begin{tabular}{lrrrr}
    \toprule
        Dataset & Train size & Test size & Length & Classes\\
        \midrule
        Beef & 30 & 30 & 470 & 5  \\
        CBF & 30 & 900 & 128 & 3  \\ 
        ChlorineConcentration & 467 & 3840 & 166 & 3\\
        CinCECGTorso & 40 & 1380 & 1639 & 4  \\ 
        Coffee & 28 & 28 & 286 & 2  \\ 
        DiatomSizeReduction & 16 & 306 & 345 & 4  \\ 
        ECG200 & 100 & 100 & 96 & 2  \\ 
        ECGFiveDays & 23 & 861 & 136 & 2  \\ 
        FaceFour & 24 & 88 & 350 & 4  \\ 
        GunPoint & 50 & 150 & 150 & 2  \\ 
        Haptics & 155 & 308 & 1092 & 5  \\ 
        ItalyPowerDemand & 67 & 1029 & 24 & 2  \\
        Lightning2 & 60 & 61 & 637 & 2\\
        Lightning7 & 70 & 73 & 319 & 7\\
        MoteStrain & 20 & 1252 & 84 & 2  \\ 
        OliveOil & 30 & 30 & 570 & 4  \\ 
        SonyAIBORobotSurface1 & 20 & 601 & 70 & 2  \\ 
        SonyAIBORobotSurface2 & 27 & 953 & 65 & 2  \\ 
        Symbols & 25 & 995 & 398 & 6  \\ 
        SyntheticControl & 300 & 300 & 60 & 6\\
        Trace & 100 & 100 & 275 & 4  \\ 
        TwoLeadECG & 23 & 1139 & 82 & 2  \\
        TwoPatterns & 1000 & 4000 & 128 & 4  \\ 
        UWaveGestureX & 896 & 3582 & 315 & 8  \\ 
        UWaveGestureY & 896 & 3582 & 315 & 8  \\ 
        UWaveGestureZ & 896 & 3582 & 315 & 8  \\ 
        Wafer & 1000 & 6164 & 152 & 2  \\ \bottomrule
    \end{tabular}
    \label{timeseries}
\end{table}

\begin{table}[htb]
  \caption{Out-of-sample accuracy computed on benchmark datasets for time series classification.}
  \label{tab:timeseries}
  \centering
  \begin{tabular}{l@{\hspace{0.5em}}*{2}{c}}
    \toprule
    Dataset & Our method & gRSF (Full model) \\
    \midrule
    Beef & 0.45~$\pm$ 0.10 & 0.70~$\pm$ 0.02 \\
    CBF & 0.83~$\pm$ 0.06 & 0.99~$\pm$ 0.00 \\
    ChlorineConcentration & 0.56 $\pm$ 0.01 & 0.58 $\pm$ 0.00 \\
    CinCECGTorso & 0.62~$\pm$ 0.05 & 0.79~$\pm$ 0.00 \\
    Coffee & 0.87~$\pm$ 0.06 & 0.90~$\pm$ 0.02 \\
    DiatomSizeReduction & 0.83~$\pm$ 0.05 & 0.96~$\pm$ 0.01 \\
    ECG200 & 0.81~$\pm$ 0.05 & 0.84~$\pm$ 0.01 \\
    ECGFiveDays & 0.97~$\pm$ 0.03 & 1.00~$\pm$ 0.00 \\
    FaceFour & 0.74~$\pm$ 0.07 & 0.98~$\pm$ 0.00 \\
    GunPoint & 0.88~$\pm$ 0.06 & 0.99~$\pm$ 0.00 \\
    Haptics & 0.39~$\pm$ 0.03 & 0.43~$\pm$ 0.02 \\
    ItalyPowerDemand & 0.93~$\pm$ 0.01 & 0.94~$\pm$ 0.00 \\
    Lightning2 & 0.77 $\pm$ 0.03 & 0.79 $\pm$ 0.01\\
    Lightning7 & 0.55 $\pm$ 0.05 & 0.66 $\pm$ 0.01\\
    MoteStrain & 0.73~$\pm$ 0.03 & 0.87~$\pm$ 0.01 \\
    OliveOil & 0.69~$\pm$ 0.05 & 0.86~$\pm$ 0.01 \\
    SonyAIBORobotSurface1 & 0.91~$\pm$ 0.05 & 0.96~$\pm$ 0.00 \\
    SonyAIBORobotSurface2 & 0.86~$\pm$ 0.04 & 0.94~$\pm$ 0.00 \\
    Symbols & 0.61~$\pm$ 0.10 & 0.92~$\pm$ 0.01 \\
    SyntheticControl & 0.93 $\pm$ 0.02 & 0.99 $\pm$ 0.00\\
    Trace & 0.95~$\pm$ 0.02 & 0.99~$\pm$ 0.00 \\
    TwoLeadECG & 0.91~$\pm$ 0.03 & 0.95~$\pm$ 0.00 \\
    TwoPatterns & 0.62~$\pm$ 0.04 & 0.86~$\pm$ 0.01 \\
    UWaveGestureX & 0.63~$\pm$ 0.03 & 0.70~$\pm$ 0.00 \\
    UWaveGestureY & 0.53~$\pm$ 0.03 & 0.62~$\pm$ 0.00 \\
    UWaveGestureZ & 0.78~$\pm$ 0.02 & 0.64~$\pm$ 0.01 \\
    Wafer & 0.96~$\pm$ 0.01 & 0.97~$\pm$ 0.00 \\
    \midrule
    \textit{Average} & 0.74~$\pm$ 0.04 & \textbf{0.84}~$\pm$ 0.01\\
    \bottomrule
  \end{tabular}
\end{table}

\normalsize

\clearpage

\section{Conclusions and future work} \label{sec6}
In this paper, we introduce a novel method to extract interpretable rules from tree ensembles. Our approach involves solving a set partitioning problem, separating the computation of loss and regularization functions from the Integer Programming model. Leveraging the flexibility of Mathematical Programming, our method can incorporate various non-linear functions, and it can handle both multi-class classification and regression tasks for tabular and temporal data.
We present extensive and statistically significant computational results for the classification and regression of tabular data, showing that the proposed method achieves comparable performance to other competing rule extraction methods for tree ensembles. \lb{We also conduct a comparison in terms of explanation fidelity, introducing two novel measures of internal fidelity for surrogate models from tree ensembles. The results demonstrate that the proposed approach outperforms others in internal fidelity, while maintaining competitive external fidelity.}
Additionally, we provide empirical evidence that the proposed method can effectively extract interpretable rules from tree ensembles designed for the multi-class classification of temporal data.
Experiments also show that our method is computationally demanding and less efficient than competing approaches. Nevertheless, this behavior is almost completely ascribable to the preprocessing steps, empirically proving the efficiency of the proposed IP model.
Future research may focus on improving the preprocessing phase and extending the implementation of the proposed method to handle tree ensembles for more complex tasks, such as image~\cite{bosch2007image, man2018image, xu2012improved} or text classification~\cite{bouaziz2014short, sun2020application}.


\section*{Acknowledgments}
This work was supported by Fedegari Autoclavi S.p.A\@. The authors thank Prof. Stefano Gualandi and Prof. Isak Samsten for their valuable support. The authors express their gratitude to all the researchers who generously provided the datasets used in this work.

\lb{\appendix
\section{Comparison of rule extraction methods}\label{appendix:A}
We note that competing methods can generate significantly different types of rules. RuleFit is the only method that produces rules similar in form to those generated by our approach. In contrast, SIRUS generates rules in an ``if-then-else'' format, where each rule could be interpreted as two separate rules according to our definition, due to the inclusion of the ``else'' part. Finally, NodeHarvest generates rules in the form of weighted trees, where each node contributes to the prediction and is weighted depending on the tree structure. In the case of NodeHarvest, branch nodes of weighted trees can be viewed as rules, as they may become leaf nodes after the pruning step. As a result, ensuring a fair comparison across methods is not trivial. Nonetheless, each of the methods considered allows for direct control over the upper limit of rules to be extracted. Therefore, in our experimental evaluation we choose to maintain the same fixed number of rules according to the corresponding definitions used by the competing methods, as specified by their authors.}

\lb{\section{Interpretable output}\label{appendix:B}
We present examples of the interpretable rules extracted by different methods. The target datasets selected for this qualitative analysis are those from the benchmark where all methods performed similarly. For simplicity, we consider RuleFit and SIRUS as competing methods, since their lists of rules are easier to represent.
\paragraph{Tecator} For this tabular dataset, we trained a random forest with maximal depth 2 and we limited the number of extracted rules up to 4. We do not apply any train-test split, as we are not evaluating out-of-sample performance.
For completeness, we report the prediction values of SIRUS rules obtained by the ``else'' clause. In Table \ref{app:tecator}, we present the lists of rules obtained from our method, RuleFit, and SIRUS, along with their corresponding weights and prediction values. Notably, the rules extracted by our method are the only ones that construct a complete decision tree, as depicted in Figure \ref{app:teca1}. Instead, RuleFit tends to extract logical conditions with slightly different thresholds, and SIRUS tends to produce shorter rules, due to the inclusion of ``else'' conditions.}

\begin{table}[htb]
\centering
\begin{tabular}{lrrr}
\toprule
\multicolumn{4}{c}{\textbf{Our method}} \\
Rule & Weight & Prediction & \\ \hline
$(\bm{x}^{(122)}, 59.20, \leq), (\bm{x}^{(122)}, 47.25, \leq)$ &  - & 2.12 &\\
$(\bm{x}^{(122)}, 59.20, \leq), (\bm{x}^{(122)}, 47.25, >)$ & - & 0.84 & \\
$(\bm{x}^{(122)}, 59.20, >), (\bm{x}^{(122)}, 67.20, \leq)$ & - & -0.13 & \\
$(\bm{x}^{(122)}, 59.20, >), (\bm{x}^{(122)}, 67.20, >)$ & - & -0.80 & \\
\midrule
\multicolumn{4}{c}{\textbf{RuleFit}} \\
Rule & Weight & Prediction  \\ \hline
$(\bm{x}^{(122)}, 57.95, \leq), (\bm{x}^{(122)}, 46.25, \leq)$ & 1.78 & 2.21 & \\
$(\bm{x}^{(122)}, 59.20, \leq), (\bm{x}^{(122)}, 47.25, >)$ & 0.54 & 0.87 & \\
$(\bm{x}^{(122)}, 57.95, >)$ & -0.45 & -0.54 & \\
$(\bm{x}^{(122)}, 59.20, >), (\bm{x}^{(122)}, 67.20, >)$ & -0.67 & -0.80 & \\ \midrule
\multicolumn{4}{c}{\textbf{SIRUS}} \\
Rule & Weight & Prediction & Prediction (\textit{else}) \\ \hline
$(\bm{x}^{(122)}, 57.97, \leq)$ & 0.57 & 1.32 & -0.57 \\
$(\bm{x}^{(122)}, 44.81, \leq)$ & 0.44 & 2.16 & -0.24\\
$(\bm{x}^{(123)}, 16.00, \leq)$ & 0.27 & 1.19 & -0.50 \\
\bottomrule
\end{tabular}
\caption{List of rules obtained by competing methods on the \textit{tecator} dataset.}
\label{app:tecator}
\end{table}

\begin{figure}[!htb]
\centering
\begin{tikzpicture}[level distance=2cm,
level 1/.style={sibling distance=7.5cm},
level 2/.style={sibling distance=3cm}]
\node[circle,draw,fill=lightgray!20,text=black] (1) {1}
child {node[circle,draw,fill=lightgray!20,text=black] (2) {2}
child {node[circle,draw,fill=lightgray!20,text=black] (4) {2.12}}
child {node[circle,draw,fill=lightgray!20,text=black] (5) {0.84}}
}
child {node[circle,draw,fill=lightgray!20,text=black] (3) {3}
child {node[circle,draw,fill=lightgray!20,text=black] (6) {-0.13}}
child {node[circle,draw,fill=lightgray!20,text=black] (7) {-0.80}}
};
\path (1) -- (2) node [pos=0.1, align=center, left=0.4cm, font=\normalsize, color=black] {$\bm{x}^{(122)} \leq 59.20$};
\path (1) -- (3) node [pos=0.1, align=center, right=0.4cm, font=\normalsize, color=black] {$\bm{x}^{(122)} > 59.20$};
\path (2) -- (4) node [pos=0.1, align=center, left=0.3cm, font=\normalsize, color=black] {$\bm{x}^{(122)} \leq 47.25$};
\path (2) -- (5) node [pos=0.1, align=center, right=0.3cm, font=\normalsize, color=black] {$\bm{x}^{(122)} > 47.25$};
\path (3) -- (6) node [pos=0.1, align=center, left=0.3cm, font=\normalsize, color=black] {$\bm{x}^{(122)} \leq 67.20$};
\path (3) -- (7) node [pos=0.1, align=center, right=0.3cm, font=\normalsize, color=black] {$\bm{x}^{(122)} > 67.20$};
\end{tikzpicture}
\caption{Rules extracted from \textit{tecator} by our method in the form of a decision tree.}
\label{app:teca1}
\end{figure}

\lb{\paragraph{Houses} For this tabular dataset, we trained a random forest with maximal depth 3 and we limited the number of extracted rules up to 15. We do not apply any train-test split, as we are not evaluating out-of-sample performance. For completeness, we report the prediction values of SIRUS rules obtained by the ``else'' clause. In Table \ref{app:houses}, we present the lists of rules obtained from our method, RuleFit, and SIRUS, along with their corresponding weights and prediction values. Again, the rules extracted by our method are the only ones that form a complete decision tree, as shown in Figure \ref{app:houses1}. Instead, RuleFit tends to extract shorter rules containing logical conditions with slightly different thresholds, while SIRUS tends to produce more rules made of only one condition.}

\begin{table}[ht]
\centering
\begin{tabular}{lrrr}
\toprule
\multicolumn{4}{c}{\textbf{Our method}} \\
Rule & Weight & Prediction & \\ \hline
$(\bm{x}^{(0)}, 5.08, \leq), (\bm{x}^{(0)}, 3.11, \leq), (\bm{x}^{(6)}, 34.45, \leq)$ & - & -0.42 & \\
$(\bm{x}^{(0)}, 5.08, \leq), (\bm{x}^{(0)}, 3.11, \leq), (\bm{x}^{(6)}, 34.45, >)$ & - & -0.79 &   \\
$(\bm{x}^{(0)}, 5.08, \leq), (\bm{x}^{(0)}, 3.11, >), (\bm{x}^{(1)}, 38.50, \leq)$ & - & -0.08 &  \\
$(\bm{x}^{(0)}, 5.08, \leq), (\bm{x}^{(0)}, 3.11, >), (\bm{x}^{(1)}, 38.50, >)$ & - & 0.45 &   \\
$(\bm{x}^{(0)}, 5.08, >), (\bm{x}^{(0)}, 6.87, \leq), (\bm{x}^{(1)}, 36.50, \leq)$ & - & 0.64 &  \\
$(\bm{x}^{(0)}, 5.08, >), (\bm{x}^{(0)}, 6.87, \leq), (\bm{x}^{(1)}, 36.50, >)$ & - & 1.35 & \\
$(\bm{x}^{(0)}, 5.08, >), (\bm{x}^{(0)}, 6.87, >), (\bm{x}^{(0)}, 7.82, \leq)$ & - & 1.45 &  \\
$(\bm{x}^{(0)}, 5.08, >), (\bm{x}^{(0)}, 6.87, >), (\bm{x}^{(0)}, 7.82, >)$ & - & 2.17 & \\
\midrule
\multicolumn{4}{c}{\textbf{RuleFit}} \\
Rule & Weight & Prediction &  \\ \hline
$(\bm{x}^{(0)}, 5.04, \leq), (\bm{x}^{(0)}, 2.85, \leq)$ & -0.33 & -0.67 & \\
$(\bm{x}^{(0)}, 5.04, >)$ & 0.38 & 1.06 & \\
$(\bm{x}^{(0)}, 5.08, \leq), (\bm{x}^{(0)}, 3.12, \leq)$ & -0.36 & -0.60 &  \\
$(\bm{x}^{(0)}, 5.08, >), (\bm{x}^{(0)}, 6.82, >)$ & 1.02 & 1.86 & \\
$(\bm{x}^{(0)}, 5.09, \leq)$ & -0.05 & -0.28 &  \\
$(\bm{x}^{(0)}, 5.37, \leq)$ & -0.34 & -0.24 & \\
$(\bm{x}^{(0)}, 5.12, \leq)$ & -0.04 & -0.27 & \\ \midrule
\multicolumn{4}{c}{\textbf{SIRUS}} \\
Rule & Weight & Prediction & Prediction (\textit{else}) \\ \hline
$(\bm{x}^{(0)}, 5.11, \leq)$ & 0.41 & -0.28 & 1.11 \\
$(\bm{x}^{(6)}, 38.5, \leq)$ & 0.45 & 0.08 & -0.75 \\
$(\bm{x}^{(7)}, 122, \leq)$ & 0.17 & 0.33 & -0.08 \\
$(\bm{x}^{(1)}, 46, \leq)$ & 0.12 & -0.04 & 0.35 \\
$(\bm{x}^{(0)}, 5.11, \leq), (\bm{x}^{(0)}, 3.14, \leq)$ & 0.53 & -0.60 & 0.40 \\
$(\bm{x}^{(2)}, 2460, \leq)$ & 0.11 & -0.13 & 0.20 \\
$(\bm{x}^{(0)}, 5.11, \leq), (\bm{x}^{(1)}, 46, \leq)$ & 0.29 & -0.33 & 0.79 \\
$(\bm{x}^{(2)}, 1840, \leq)$ & 0.12 & -0.20 & 0.13 \\
$(\bm{x}^{(2)}, 2130, \leq)$ & 0.11 & -0.16 & 0.16 \\
$(\bm{x}^{(0)}, 5.11, \geq), (\bm{x}^{(0)}, 3.14, \leq), (\bm{x}^{(6)}, 38.5, \leq)$ & 0.07 & 0.10 & -0.06 \\
$(\bm{x}^{(0)}, 5.11, \leq), (\bm{x}^{(7)}, 122, \geq)$ & 0.21 & -0.34 & 0.63 \\
\bottomrule
\end{tabular}
\caption{List of rules obtained by competing methods on the \textit{houses} dataset.}
\label{app:houses}
\end{table}

\begin{figure}[htb]
\centering
\normalsize
\begin{tikzpicture}[rotate=90, scale=0.85, transform shape, level distance=2cm,
level 1/.style={sibling distance=12cm},
level 2/.style={sibling distance=6cm},
level 3/.style={sibling distance=2cm}]
\node[circle,draw,fill=lightgray!20,text=black] (1) {1}
child {node[circle,draw,fill=lightgray!20,text=black] (2) {2}
child {node[circle,draw,fill=lightgray!20,text=black] (4) {3}
child {node[circle,draw,fill=lightgray!20,text=black] (8) {-0.42}}
child {node[circle,draw,fill=lightgray!20,text=black] (9) {-0.79}}}
child {node[circle,draw,fill=lightgray!20,text=black] (5) {4}
child {node[circle,draw,fill=lightgray!20,text=black] (10) {-0.08}}
child {node[circle,draw,fill=lightgray!20,text=black] (11) {0.45}}}
}
child {node[circle,draw,fill=lightgray!20,text=black] (3) {5}
child {node[circle,draw,fill=lightgray!20,text=black] (6) {6}
child {node[circle,draw,fill=lightgray!20,text=black] (12) {0.64}}
child {node[circle,draw,fill=lightgray!20,text=black] (13) {1.35}}}
child {node[circle,draw,fill=lightgray!20,text=black] (7) {7}
child {node[circle,draw,fill=lightgray!20,text=black] (14) {1.45}}
child {node[circle,draw,fill=lightgray!20,text=black] (15) {2.17}}}
};

\path (1) -- (2) node [pos=-0.1, align=center, left=0.6cm, font=\small, color=black] {$\bm{x}^{(0)} \leq 5.08$};
\path (1) -- (3) node [pos=-0.1, align=center, right=0.6cm, font=\small, color=black] {$\bm{x}^{(0)} > 5.08$};

\path (2) -- (4) node [pos=0.1, align=center, left=0.3cm, font=\small, color=black] {$\bm{x}^{(0)} \leq 3.11$};
\path (2) -- (5) node [pos=0.1, align=center, right=0.3cm, font=\small, color=black] {$\bm{x}^{(0)} > 3.11$};
\path (3) -- (6) node [pos=0.1, align=center, left=0.3cm, font=\small, color=black] {$\bm{x}^{(0)} \leq 6.87$};
\path (3) -- (7) node [pos=0.1, align=center, right=0.3cm, font=\small, color=black] {$\bm{x}^{(0)} > 6.87$};

\path (4) -- (8) node [pos=0.1, align=center, left=0.3cm, font=\small, color=black] {$\bm{x}^{(6)} \leq 34.45$};
\path (4) -- (9) node [pos=0.1, align=center, right=0.3cm, font=\small, color=black] {$\bm{x}^{(6)} > 34.45$};
\path (5) -- (10) node [pos=0.1, align=center, left=0.3cm, font=\small, color=black] {$\bm{x}^{(1)} \leq 38.5$};
\path (5) -- (11) node [pos=0.1, align=center, right=0.3cm, font=\small, color=black] {$\bm{x}^{(1)} > 38.5$};
\path (6) -- (12) node [pos=0.1, align=center, left=0.3cm, font=\small, color=black] {$\bm{x}^{(1)} \leq 36.5$};
\path (6) -- (13) node [pos=0.1, align=center, right=0.3cm, font=\small, color=black] {$\bm{x}^{(1)} > 36.5$};
\path (7) -- (14) node [pos=0.1, align=center, left=0.3cm, font=\small, color=black] {$\bm{x}^{(0)} \leq 7.82$};
\path (7) -- (15) node [pos=0.1, align=center, right=0.3cm, font=\small, color=black] {$\bm{x}^{(0)} > 7.82$};

\end{tikzpicture}
\caption{Rules extracted from \textit{house} by our method in the form of a decision tree.}
\label{app:houses1}
\end{figure}

\lb{\paragraph{ECGFiveDays}
This temporal dataset contains ECG signals recorded from a 67-year-old male. The two classes correspond to the dates 12/11/1990 and 17/11/1990 \cite{UCRArchive2018}. For this dataset, we first trained the gRSF over the training set by using the provided hyperparameters \cite{KarlssonForest} and setting the maximal depth to 3. Since no competing methods currently address time series classification, we applied our rule extraction method with the $\ell$ parameter fixed at 2. The resulting list of rules is shown in Figure \ref{app:ecg}, presented as a shapelet-based decision tree with unitary depth. The first leaf node corresponds to the date 17/11/1990, while the second leaf node corresponds to 12/11/1990. We also present two examples of time series from different classes, overlaying the single shapelet extracted from the gRSF and reporting the distance values computed using the subsequence distance measure \eqref{eq:subdist}. Notably, the temporal positions of the two peaks near $x=80$ are sufficient for accurately categorizing the signals. The extracted rules effectively capture and explain the predictions of the shapelet-based tree ensemble, obtaining a test accuracy of $99\%$.}

\begin{figure}[htb]
\centering

\begin{tikzpicture}[level distance=2.5cm,
level 1/.style={sibling distance=4.5cm}]

\node[circle,draw,fill=lightgray!20,text=black] (1) {1}
child {node[circle,draw,fill=lightgray!20,text=black] (2) {17/11}}
child {node[circle,draw,fill=lightgray!20,text=black] (3) {12/11}};

\path (1) -- (2) node [pos=0.1, align=center, left=0.4cm, font=\normalsize, text=black] {dist$(\bm{x}, \bm{s}_1) \leq 3.53$};
\path (1) -- (3) node [pos=0.1, align=center, right=0.4cm, font=\normalsize, text=black] {dist$(\bm{x}, \bm{s}_1) > 3.53$};

\end{tikzpicture}

\vspace{1cm} 

\begin{tikzpicture}

\begin{axis}[
    width=7cm, height=4.5cm, 
    xlabel={dist$(\bm{x}_3, \bm{s}_1) = 2.60$},
    axis x line=bottom,
    axis y line=left,
    ytick={-6, -4, -2, 0, 2}, 
    ymin=-6, ymax=3, 
    legend style={at={(1,0.10)}, anchor=south east, legend columns=1},
    every axis plot/.append style={thick},
]

\addplot[red, ultra thick, mark=none] table {shapelet_2.txt};
\addlegendentry{$\bm{s}_1$}
\addplot[blue, mark=none] table {timeseries_2.txt};
\addlegendentry{$\bm{x}_3$}
\end{axis}

\begin{scope}[xshift=7cm]  
\begin{axis}[
    width=7cm, height=4.5cm, 
    xlabel={dist$(\bm{x}_4, \bm{s}_1) = 4.85$},
    axis x line=bottom,
    axis y line=left,
    ytick={-6, -4, -2, 0, 2}, 
    ymin=-6, ymax=3, 
    legend style={at={(1,0.10)}, anchor=south east, legend columns=1},
    every axis plot/.append style={thick},
]

\addplot[red, ultra thick, mark=none] table {shapelet.txt};
\addlegendentry{$\bm{s}_1$}

\addplot[blue, mark=none] table {timeseries.txt};
\addlegendentry{$\bm{x}_4$}
\end{axis}
\end{scope}

\end{tikzpicture}

\caption{Rules extracted from \textit{ECGFiveDays} by our method in the form of a shapelet-based decision tree.}
\label{app:ecg}
\end{figure}


\clearpage







\bibliographystyle{cas-model2-names}

\bibliography{main}

\end{document}